\documentclass[12pt,leqno]{amsart}
\usepackage{amsmath}
\usepackage{amssymb} 
\usepackage{mathrsfs}
\usepackage[all]{xy}
\usepackage{color}
\usepackage{soul}
\usepackage[hidelinks]{hyperref}
\usepackage{enumitem}
\usepackage{graphicx}

\usepackage{physics}
\usepackage{amsmath}
\usepackage{tikz}
\usepackage{mathdots}
\usepackage{yhmath}
\usepackage{cancel}
\usepackage{color}
\usepackage{siunitx}
\usepackage{array}
\usepackage{multirow}
\usepackage{amssymb}
\usepackage{gensymb}
\usepackage{tabularx}
\usepackage{extarrows}
\usepackage{booktabs}
\usetikzlibrary{fadings}
\usetikzlibrary{patterns}
\usetikzlibrary{shadows.blur}
\usetikzlibrary{shapes}

\newcommand\N{{\mathbb N}}
\newcommand\R{{\mathbb R}}
\newcommand\T{{\mathbb T}}

\newcommand\Z{{\mathbb Z}}

\newcommand\Sp{{\mathbb S}}

\def\AA{{\mathcal A}}
\def\BB{{\mathcal B}}
\def\CC{{\mathcal C}}

\def\EE{{\mathcal E}}

\def\LL{{\mathcal L}}

\def\PP{{\mathcal P}}

\def\RR{{\mathcal R}}
\def\SS{{\mathcal S}}
\def\TT{{\mathcal T}}

\def\XX{{\mathcal X}}

\def\BBB{{\mathscr B}}
\def\CCC{{\mathscr C}}

\def\XXX{{\mathscr X}}

\def\eps{{\varepsilon}}

\newtheorem{theo}{Theorem}
\newtheorem{prop}{Proposition}
\newtheorem{lem}[prop]{Lemma}
\newtheorem{cor}[prop]{Corollary}
\newtheorem{rem}[prop]{Remark}

\newtheorem{notation}[prop]{Notations}

\newcommand{\resetcounter}{
	\setcounter{equation}{0}
	\setcounter{theo}{0}
	\setcounter{prop}{0}
}

\newcommand{\step}[2]{\medskip\noindent\textit{Step #1: #2.}}

\newcommand{\hyproj}{\mathbf{P}^\eps_\flat}
\newcommand{\hyprojo}{\mathbf{P}^\eps_\sharp}

\newcommand{\vertiii}[1]{{\left\vert\kern-0.25ex\left\vert\kern-0.25ex\left\vert #1 
		\right\vert\kern-0.25ex\right\vert\kern-0.25ex\right\vert}}

\newcommand{\nul}{\textnormal{N}}

\newcommand{\Id}{\textnormal{Id} \,}

\newcommand{\WP}{\textnormal{WP}}
\newcommand{\IP}{\textnormal{IP}}
\newcommand{\disp}{\textnormal{disp}}
\newcommand{\ac}{\textnormal{ac}}

\newcommand{\ini}{\textnormal{in}}
\newcommand{\solp}{\overline{h}}
\newcommand{\solg}{\overline{g}}

\usepackage{etoolbox}

\newcommand{\spcp}{\EE}
\newcommand{\Spcp}[2]{{\EE^{#1, #2}}}
\newcommand{\spcpp}{{\EE_\nu}}
\newcommand{\Spcpp}[2]{{\EE_\nu^{{#1, #2}} } }
\newcommand{\tspcp}{\XX_\eps}
\newcommand{\tSpcp}[2]{ \XX^{#1,#2, \sigma}_{\eps} }
\newcommand{\tspcpN}[1]{\left\|#1\right\|_{\tspcp }}
\newcommand{\tSpcpN}[3]{\left\|#1\right\|_{\tSpcp{#2}{#3} }}

\newcommand{\spcg}{E}
\newcommand{\Spcg}[1]{{\spcg^{#1}} }
\newcommand{\tspcg}[2]{ {X_{#1}^{#2}} }
\newcommand{\tspcgN}[3]{\left\|#1\right\|_{\tspcg{#2}{#3}} }
\newcommand{\spcgg}{H} 

\renewcommand{\d}{\textnormal{d}}

\newcommand{\err}{\mathbf{r}}
\newcommand{\ah}{\mathbf{a}}

\title[From Boltzmann to Navier-Stokes-Fourier with general initial data]{On the convergence from Boltzmann to Navier-Stokes-Fourier for general initial data}

\begin{document}

	\begin{abstract}
		In this work, we prove the convergence of strong solutions of the Boltzman equation, for initial data having polynomial decay in the velocity variable, towards those of the incompressible Navier-Stokes-Fourier system. We show in particular that the solutions of the rescaled Boltzmann equation do not blow up before their hydrodynamic limit does. This is made possible by adapting the strategy from \cite{BMM} of writing the solution to the Boltzmann equation as the sum a part with polynomial decay and a second one with Gaussian decay. The Gaussian part is treated with an approach reminiscent of the one from~\cite{IGT20}.
	\end{abstract}
	
	\author{{\sc Pierre Gervais}}
	
	\date\today
	
	\maketitle
	
	\vspace{1cm}

	\tableofcontents
	
	\bigskip
	
	\section{Introduction}
	
	\resetcounter
	
	\subsection{The models}
	Let us present the main models considered in this paper, namely the Boltzmann equation and the incompressible Navier-Stokes-Fourier system, and how they are related.
	
	\subsubsection{The Boltzmann equation}
	\label{scn:boltzmann}
	
	Consider a rarefied gas whose particles number density at position $x~\in~\Omega$ (here $\Omega= \R^d$ or $\T^d$, and $d = 2, 3$), traveling at velocity~$v \in \R^d$ at time $t \in \R^+$ is given by~$F(t, x, v)$. Assuming the particles undergo hard spheres collisions, $F$ evolves according to the \textit{Boltzmann equation}
	\begin{equation}
		\label{eqn:boltz}
		\begin{cases}
					\left(\partial_t + v \cdot \nabla_x\right)F = Q(F,F),\\
					F_{| t = 0} = F_\ini,
		\end{cases}
	\end{equation}
	which is a transport equation whose source term models the effect of binary collisions between pairs of particles. The operator $Q$ is called the \textit{Boltzmann operator} or \textit{collision operator} and is an integral bilinear symmetric operator defined as
	\begin{equation*}
		Q(F, G) := \frac{1}{2} (Q^+(F, G) - Q^-(F, G) + Q^+(G, F) - Q^-(G, F)),
	\end{equation*}
	where the so-called \textit{gain} and \textit{loss part} are defined as
	\begin{gather*}
		Q^+(F,G)(v) := \int_{\R^d} \int_{\Sp^{d-1}_\sigma} |v-v_*|F' G'_* \d v_* \d\sigma,\\
		Q^-(F,G)(v) := \int_{\R^d} \int_{\Sp^{d-1}_\sigma} |v-v_*|
		F G_* \d v_* \d\sigma,
	\end{gather*}
	where we used the standard notation
	\begin{itemize}
		\item $v$ and $v_*$ for the velocities of two particles after the collision,
		\item $v'$ and $v_*'$ for their velocities before the collision, given by
		\begin{equation}
			\label{eqn:pre_collision_velocities}
			v' = \frac{v + v_*}{2} + \frac{|v-v_*|}{2} \sigma, ~ v_*' = \frac{v + v_*}{2} - \frac{|v-v_*|}{2} \sigma,
		\end{equation}
		\item $F:=F(v)$, $F':=F(v')$, $G_*':=G(v_*')$ and $G_*:=G(v_*)$.
	\end{itemize}
	One can retrieve local macroscopic observables of the gas such as the mass density~$R(t, x)$, local temperature $T(t, x)$ and bulk velocity $U(t, x) \in \R^d$ using the moments of the particles number density:
	\begin{gather*}
		R(t, x) := \int_{\R^d} F(t, x, v) \d v,\\
		(RU)(t, x) := \int_{\R^d} v F(t, x, v) \d v,\\
		R(|U|^2 + d T)(t, x) := \int_{\R^d} |v|^2F(t, x, v) \d v.
	\end{gather*}
	Using the symmetries of $Q$ which imply that $Q(F, F)$ is orthogonal to $1, v, |v|^2$, one may show by integrating \eqref{eqn:boltz} against $1, v, |v|^2$ in $v$ that these quantities satisfy the following conservation laws:
	\begin{gather*}
		\partial_t R(t, x) + \nabla_x \cdot U(t, x) = 0,\\
		\partial_t U(t, x) + \nabla_x \cdot \int_{\R^d} v \otimes v F(t, x, v) \d v = 0,\\
		\partial_t (R |U|^2 + d R T)(t, x) + \nabla_x \cdot \int_{\R^d} v |v|^2 F(t, x, v) \d v = 0.
	\end{gather*}
	The last physical quantity we wish to introduce is the \textit{entropy} of the gas, defined as
	\begin{equation*}
		H(t) := \int_{\R^d \times \Omega} F(t, x, v) \log F(t, x, v) \d v \d x.
	\end{equation*}
	The celebrated \textit{H-theorem} (\textit{second law of theormodynamics}) states that~$H$ is a non-negative non-increasing function and thus a Lyapunov functionnal. The entropy minimizers are therefore equilibria, called \textit{Maxwellian distributions} (see for instance \cite[Chapter 1]{UY06}), which are Gaussians in the variable $v$, parametrized by the macroscopic observables:
	\begin{equation*}
		M_{R, T, U}(t, x, v) := \frac{R(t, x)}{(2 \pi T(t, x))^{d/2}} \exp \left(-\frac{|v - U(t, x)|^2}{2 T(t, x)}\right).
	\end{equation*}
	In this paper, we will denote the normal centered distribution ($R=1, U=0, T=1$) by
	$$M(v) := \frac{1}{\left(2 \pi\right)^{d/2}} \exp\left(- \frac{|v|^2}{2}\right),$$
	and we define the corresponding so-called \textit{collision frequency}
	\begin{equation}
		\label{eqn:def_coll_freq}
		\nu(v) := \int_{\R^d} M_* |v-v_*| \d v_*,
	\end{equation}
	which is such that for some constants $0 < \nu_0 \leq \nu_1$, denoting $\langle v \rangle := \sqrt{1 + |v|^2}$,
	\begin{equation}
		\label{eqn:bound_col_freq}
		\nu_0 \langle v \rangle \leq \nu(v) \leq \nu_1 \langle v \rangle.
	\end{equation}
	Global weak solutions to \eqref{eqn:boltz}, called \textit{renormalized solutions}, were constructed by R. DiPerna and P. L. Lions in \cite{DL89} under very general assumptions, corresponding to initial data with finite mass, energy and entropy.
	
	In \cite{U86}, strong solutions to \eqref{eqn:boltz} have been constructed for smooth (in $x$) initial data with Gaussian decay. These solutions are global when the initial data is close to some absolute Maxwellian $M_{R, U, T}$ ($R, T > 0$, $U \in \R^d$), and local in time otherwise. In \cite{GMM17}, the authors exhibited strong solutions to \eqref{eqn:boltz} for smooth (in $x$) initial data having algebraic decay and close to some absolute Maxwellian. Let us state their existence theorem more precisely.
	\begin{theo}
		\label{thm:strong_solutions_polynomial}
		Let $\alpha > 2$, there exists some positive $\eta=\eta(\alpha)$ such that, for any initial data $F_\ini=F_\ini(x, v)$ with $x \in \T^3$ and $v \in \R^3$ satisfying
		$$\|F_\ini - M\|_{L^1_v L^\infty_x \left( \langle v \rangle^\alpha \d v \right)} < \eta(\alpha),$$
		there exists a unique global solution $F \in \CC^0_t L^1_v L^\infty_x\left( \langle v \rangle^\alpha \d v \right) \cap L^1_{t, v} L^\infty_x\left( \langle v \rangle^{\alpha+1} \d v \right)$ to the Boltzmann equation~\eqref{eqn:boltz}.
	\end{theo}
	We also mention \cite{G04} in which perturbative solutions are constructed using a non-linear energy method based on the so-called micro-macro decomposition (see Section \ref{scn:micro_macro_decomposition}), and \cite{IS84} in which solutions are constructed near vacuum.

	\subsubsection{Micro-macro decomposition}
	\label{scn:micro_macro_decomposition}
	
	A gas at thermodynamic equilibrium and close to $M$, that is to say a thermodynamic equilibrium of the form $M_{1 + \rho, u, 1 + \theta}$ for some small $\rho, u, \theta$, behaves like the sum of $M$ and a so-called \textit{infinitesimal Maxwellian}:
	\begin{equation*}
		M_{1 + \rho, u, 1 + \theta}(v) = M(v) + \left(\rho + u \cdot v + \frac{1}{2} \left(|v|^2 - d\right) \theta \right) M(v) + o(\rho, u, \theta).
	\end{equation*}
	More generally, any fluctuation $F(x,v)=M(v)+f(x, v)$ around $M$ admits a unique \textit{micro-macro decomposition} (with respect to $M$)
	\begin{gather}
		\label{eqn:micro_macro_decomp}
		f(x, v) = \left(\rho_f(x) + u_f(x) \cdot v + \frac{1}{2} \left( |v|^2 - d \right) \theta_f(x) \right) M(v) + f_\bot(x, v),\\
		\notag
		\rho_f(x) := \int_{\R^d} f(x) \d v, ~ u_f(x) := \int_{\R^d} f(x, v) v \d v,\\
		\notag
		\theta_f(x) := \frac{1}{d} \int_{\R^d} f(x, v) \left( |v|^2 - d \right) \d v
	\end{gather}
	where we call \textit{the macroscopic part of $f$} the part
	\begin{gather}
		\label{eqn:def_macro}
		\left(\rho_f(x) + u_f(x) \cdot v + \frac{1}{2} (|v|^2-d) \theta_f(x)\right) M(v),
	\end{gather}
	which is given by the projection
	\begin{gather}
		\label{eqn:def_Pi}
		\Pi f(x, v) := \sum_{j=0}^{d+1} \frac{\varphi_j(v) M(v)}{\| \varphi_j M^{1/2} \|_{L^2_v}^2} \int_{\R^d} f(x, v_*) \varphi_j(v_*) \d v_*,\\
		\notag
		\varphi_0(v) = 1, ~ \varphi_j(v) = v_j, ~ \varphi_{d+1}(v) = |v|^2,
	\end{gather}
	and $f_\bot$ is the \textit{microscopic part of $f$}, thus characterized by the fact that for almost any~$x \in \Omega$,
	\begin{equation*}
		\int_{\R^d} f_\bot(x, v) \varphi_j(v) \d v = 0, ~ j = 0,\dots, d+1.
	\end{equation*}
	We say that $f$ is \textit{well-prepared} if it is macroscopic and satisfies furthermore the incompressibility condition and the Boussinesq relation of \eqref{eqn:INSF}:
	\begin{gather*}
		\nabla_x \cdot u_f = 0,\\
		\nabla_x \left(\rho_f + \theta_f\right) = 0.
	\end{gather*}
	To sum up, any $f(x, v)$ can be decomposed uniquely as
	\begin{align}
		\label{eqn:initial_data_split}
		f = \Pi f + f_\bot= f_\WP + f_\IP + f_{\bot},
	\end{align}
	where, denoting $\mathbb{P}$ the (Leray) projector on incompressible fields, the \textit{well-prepared part of $f$} writes
	\begin{gather}
		\label{eqn:well_prepared_initial data}
		f_\WP(x, v) := \left( \overline{\rho}_f(x) + \overline{u}_f(x) \cdot v + \frac{1}{2}\left( |v|^2 - d \right) \overline{\theta}_f(x) \right),\\
		\notag
		\overline{u} := \mathbb{P} u, ~ \overline{\rho}_f := - \overline{\theta}_f := \frac{2 \rho_f - d \theta_f}{2},
	\end{gather}
	and $f_\IP := \Pi f - f_{\WP}$ is called the \textit{ill-prepared part}.

	\begin{rem}
		Note that the micro-macro decomposition is orthogonal in $L^2_v\left(M^{-1} \d v\right)$.
	\end{rem}

	\begin{rem}
		Note that, in the following, the Boussinesq relation $\nabla_x (\rho + \theta) = 0$ is equivalent to $\rho + \theta = 0$ as we will assume $(\rho, \theta)$ to be mean-free when $\Omega = \T^d$ or integrable when $\Omega = \R^d$.
	\end{rem}

	\subsubsection{The hydrodynamic model}
	
	At the macroscopic level, the dynamics of an incompressible fluid is encoded at any point $x \in \Omega$ and any time $t \in \R^+$ by its bulk velocity $u(t, x) \in \R^d$, its mass density $\rho(t, x)$, its temperature $\theta(t, x)$, as well as its pressure $p(t, x)$ (which can be deduced from $u$ in the case of an incompressible fluid). These quantities evolve according to the incompressible  Navier-Stokes-Fourier system:
	\begin{equation}
		\label{eqn:INSF}
		\tag{INSF}
		\left\{
		\begin{aligned}
			\partial_t u + u \cdot \nabla u - \mu \Delta u &= - \nabla p,\\
			\partial_t \theta + u \cdot \nabla \theta - \kappa \Delta \theta &= 0,\\
			\nabla \cdot u &= 0,\\
			\nabla (\theta + \rho) &= 0,\\
			(u, \theta, \rho)_{| t = 0} &= (u_\ini, \rho_\ini, \theta_\ini).
		\end{aligned}
		\right.
	\end{equation}
	In this system, the positive scalars $\mu$ and $\kappa$ denote respectively the kinematic viscosity and heat conductivity coefficients. Global weak solutions were shown to exist by J. Leray \cite{L34} under minimal physical assumptions on the initial data, and strong solutions were conctructed by H. Fujita and T. Kato \cite{FK64} for smooth initial data. Let us present an existence result for \eqref{eqn:INSF} which is by nomean optimal but sufficient for the study led in this paper and refer to \cite{C92, LR02, LR16, FK64} for a proof.
	
	\begin{theo}
		\label{thm:existence_INSF}
		For any $s \geq d/2-1$ and $(\rho_\ini, u_\ini, \theta_\ini) \in H^s_x$ satisfying
		\begin{gather*}
			\nabla_x\cdot u_\ini = 0,\\
			\nabla_x(\rho_\ini + \theta_\ini) = 0,
		\end{gather*}
		there exists a time $T \in (0, \infty]$ such that the system \eqref{eqn:INSF} has a unique solution in $L^\infty_t \left( [0, T) ; H^s_x \right)\cap L^2_t \left( [0, T) ; H^{s+1}_x\right)$ associated with the initial data $(\rho_\ini, u_\ini, \theta_\ini)$. Furthermore, there holds $T = \infty$ if the initial data is small enough (with respect to~$\mu$).
	\end{theo}
	
	\subsubsection{Hydrodynamic limits}
	
	By choosing a system of reference values for length, time and velocity (see for example~\cite{G05, SR09}), we obtain a dimensionless version of~the~equation:
	\begin{equation*}
		\text{Ma} \, \partial_t F+v \cdot \nabla_x F = \frac{1}{\text{Kn}} Q(F, F),
	\end{equation*}
	where $\text{Kn}$ denotes the \textit{Knudsen number} (the inverse of the average number of collisions per unit of time) and $\text{Ma}$ the \textit{Mach number} (the ratio of the bulk velocity to the speed of sound, characterizing the compressibility of the fluid). To relate the Boltzmann equation to hydrodynamic models, several formal methods were proposed, first by D. Hilbert \cite{H1900}, Chapman-Enskog \cite{C52}, H. Grad \cite{G63}, then made rigorous by several works such as \cite{BU91, MEL89, C80, LN87}. At the beginning of the nineties, a systematic approach was presented in \cite{BGL91, BGL93}: C. Bardos, F. Golse and D. Levermore showed that the only possible point of accumulation for the renormalized solutions of \eqref{eqn:bolt_scl} when $\eps$ goes to zero are global weak solutions to Euler or Navier-Stokes equations. In the case of the Navier-Stokes-Equations, choosing $\text{Ma} = \text{Kn} = \eps$ in the previous scaled Boltzmann equation and performing the linearization $F=F^\eps = M + \eps f^\eps$, the authors considered the following equation:
	\begin{gather}
		\label{eqn:bolt_scl}
		\tag{B$^\eps$}
		\begin{cases}
			\partial_t f^\eps = \displaystyle \LL^\eps f^\eps + \frac{1}{\eps} Q(f^\eps, f^\eps),\\
			f^\eps(0) = f_\ini,
		\end{cases}\\
		\label{eqn:def_L_eps}
		\LL^\eps := \eps^{-2} \left(\LL - \eps v \cdot \nabla_x\right),\\
		\label{eqn:def_L}
		\LL h := 2Q(M, h),
	\end{gather}
	and showed that if the fluctation profile $f^\eps$ converges in some weak sense to some~$f^0$ as $\eps$ goes to zero, then it must write
		\begin{equation}
		\label{eqn:kin_counter}
		f^0(t, x, v) := \left(\rho^0(t, x) + u^0(t, x) \cdot v + \frac{1}{2}(|v|^2 - d) \theta^0(t, x)\right) M(v),
	\end{equation}
	and the coefficients $\rho^0, u^0$ and $\theta^0$ are distributional solutions of \eqref{eqn:INSF} with the initial condition
	\begin{gather*}
		\left(\rho^0, u^0,\theta^0 \right)_{| t = 0} = \left(\frac{2 \rho_\ini - d \theta_\ini}{d+2}, \mathbb{P} u_\ini, \frac{- 2 \rho_\ini + d \theta_\ini}{d+2}\right),\\
		\rho_\ini(x) := \int_{\R^d} f_\ini(x, v) \d v, ~
		u_\ini(x) := \int_{\R^d} v f_\ini(x, v) \d v,\\
		\theta_\ini(x) := \frac{1}{d} \int_{\R^d} (|v|^2 - d )f_\ini(x, v) \d v.
	\end{gather*}
	Recalling the definition of Section \ref{scn:micro_macro_decomposition}, this means  $f^0_{|t = 0} = f_{\ini, \WP}$.
	Finally, the diffusion coefficients in \eqref{eqn:INSF} can be expressed as
	\begin{equation*}
		\mu := \frac{1}{(d-1)(d+2)} \int_{\R^d} \text{trace} \bigg(\Phi \left(\LL \Phi\right)^T \bigg) \d v, ~ \kappa := \frac{2}{d(d+2)} \int_{\R^d} \phi \LL \phi \d v,
	\end{equation*}
	where $\Phi$ and $\phi$ are the unique functions satisfying
	\begin{gather*}
		\LL \Phi = \left(\frac{|v|^2}{2} \text{Id} - v \otimes v\right) M, ~ \LL \phi = v \left(\frac{d+2}{2} - \frac{|v|^2}{2}\right) M,
	\end{gather*}
	and orthogonal to the null-space of $\LL$ (see \cite{UY06}), which constists of macroscopic distributions:
	\begin{equation}
		\label{eqn:def_nul_L}
		\nul(\LL) := \left\{v \mapsto \left( a + b \cdot v + c (|v|^2-d)\right) M(v) , ~ a, c \in \R, \, b \in \R^d \right\}.
	\end{equation}
	This means that the heat conductivity and kinematic viscosity only depend on the physical properties of the reference equilibrium $M$ and the collision operator~$Q$. In \cite{GSR04, GSR09}, F. Golse and L. Saint-Raymond proved the compactness of these solutions and thus the convergence. We also mention \cite{LM01} in which different interactions between particles are considered.
	\medskip
	
	The derivation of the Navier-Stokes system from the Boltzmann equation in the case of strong solutions was first considered by C. Bardos and S. Ukai in \cite{BU91}, for initial data having Gaussian decay. Their approach relied on the properties of the homogeneous linearized operator $\LL$ studied by H. Grad in \cite{G63}, and on the spectral study of the full linearized operator $\LL + v \cdot \nabla_x$ led by R. Ellis and M. Pinsky in \cite{EP75}, completed by S. Ukai in \cite{U86}. They considered \eqref{eqn:bolt_scl} in its integral form
	\begin{equation}
		\label{eqn:int_bolt_scl}
		f^\eps = U^\eps f_\ini + \Psi^\eps(f^\eps, f^\eps),
	\end{equation}
	where we have denoted
	\begin{gather}
		\notag
		U^\eps(t) := \exp(t \LL^\eps),\\
		\label{eqn:def_psi_eps}
		\Psi^\eps(f^\eps, f^\eps)(t) := \frac{1}{\eps} \int_{0}^{t} U^\eps(t-t') Q\left(f^\eps(t'), f^\eps(t')\right) \d t',
	\end{gather}
	and proved continuity bounds (uniformly in $\eps$) for $U^\eps$ and $\Psi^\eps$ in some space of functions with Gaussian decay, which allowed to prove existence of strong global solutions $f^\eps$ to \eqref{eqn:int_bolt_scl} through a fixed point argument for small initial data. They also proved the convergence of~$U^\eps$ and $\Psi^\eps$ as $\eps$ goes to zero to some $U^0$ and $\Psi^0$, which implied strong convergence  convergence of $f^\eps$ to some $f^0$ satisfying
	\begin{gather}
		\label{eqn:KINSF}
		\tag{KINSF}
		f^0 = U^0 f^0_\ini + \Psi^0(f^0, f^0),
	\end{gather}
	and we recall once again that $f^0_{| t = 0} = f^0_\ini = f_{\ini, \WP}$. We could call this equation the \textit{kinetic formulation of the Navier-Stokes-Fourier system}, as the results from \cite{BGL91, BGL93} imply that~$f^0$ is the same as \eqref{eqn:kin_counter} and thus describes the unique strong solution of \eqref{eqn:INSF} (in the sense of Theorem \ref{thm:existence_INSF}).
	
	In \cite{IGT20}, I. Gallagher and I. Tristani improved this fixed point approach by considering the equation satisfies by $f^\eps - f^0$, where $f^0$ is known to exist thanks to Theorem \ref{thm:existence_INSF}. This allowed to consider large initial data $f_\ini$ for the Boltzmann equation, and thus large $(\rho_\ini^0, u_\ini^0, \theta_\ini^0)$ for the Navier-Stokes-Fourier system, and in particular, showed that the solution $f^\eps$ to \eqref{eqn:int_bolt_scl} does not blow up before $f^0$ does. We state here their result.
	\begin{theo}
		Let $s > d/2$, $\beta > d/2 + 1$, and $f_\ini \in L^\infty_v H^s_x \left( M^{-1/2} \langle v \rangle^\beta \d v \right)$. Denote~$f^0$ the solution to \eqref{eqn:KINSF} with initial data $f^0_\ini = f_{\ini, \WP}$ given by Theorem \ref{thm:existence_INSF} on a time interval~$[0, T)$. For any small enough $\eps > 0$, there exists a unique solution to \eqref{eqn:bolt_scl}
		$$f^\eps \in \CC_b \left( [0, T) ;  L^\infty_v H^s_x \left( M^{-1/2} \langle v \rangle^\beta \d v \right) \right)$$
		where $\CC_b(X;Y)$ denotes the set of bounded continuous functions from $X$ to $Y$, and it converges to $f^0$
		\begin{itemize}
			\item in $L^\infty_t L^\infty_v H^s_x \left( M^{-1/2} \langle v \rangle^\beta \d v \right)$ if $f_\ini = f_{\ini, \WP}$,
			\item in $L^\infty_t L^\infty_v H^s_x \left( M^{-1/2} \langle v \rangle^\beta \d v \right) + L^p_t L^\infty_v \left(H^s_x + W^{s, \infty} \right) \left( M^{-1/2} \langle v \rangle^\beta \d v \right)$ otherwise, where $\frac{2}{d-1} < p < \infty$.
		\end{itemize}
	\end{theo}
	
	We also mention works led in similar functional spaces: \cite{LN87, MEL89, Briant, G06, JXZ14} in which strong hydrodynamic limits of the Boltzmann equation are considered, and \cite{R20, CRT21} in which strong solutions to the incompressible Navier-Stokes-Fourier system are derived from the Landau equation.
	
	\smallskip
	The case of strong solutions to the Boltzmann equation in polynomial spaces is not different because Grad's decomposition of the linearized operator $\LL$ does not have the same nice properties as in Gaussian spaces. In particular, $\Psi^\eps$ can not be shown to be a bounded operator in $L^\infty_t L^1_v \left( \langle v \rangle^\alpha \d v \right)$ using the same approach. The study of strong solutions to the Boltzmann equation for polynomially decaying initial data was initiated by M.P. Gualdani, S. Mischler and C. Mouhot who constructed in \cite{GMM17} solutions to the Boltzmann equation for initial data close to a global Maxwellian (see Theorem \ref{thm:strong_solutions_polynomial}), or in other words (by a scaling argument), to \eqref{eqn:bolt_scl} for $f_\ini$ small and for~$\eps = 1$. This was made possible thanks to the development of the so-called \textit{enlargement theory} started by C. Mouhot in \cite{M05} and developped with S. Mischler and M.P. Gualdani in \cite{MM16, GMM17} which (loosely speaking) allows to prove spectral properties for a linear operator in large spaces under the condition that these properties hold in a smaller space and the linear operator can be decomposed in a certain~way.
	
	 In \cite{BMM}, M. Briant, S. Merino and C. Mouhot proved the solutions constructed in~\cite{GMM17} were bounded uniformly in $\eps$, and thus, up to an extraction, converging weakly to a solution of the incompressible Navier-Stokes-Fourier system.
	 
	 We also mention \cite{ALT} in which weak solutions to a modified Navier-Stokes-Fourier system are derived from solutions to the Boltzmann equation for ganular media with initial data having polynomial decay.
	 
	 \smallskip
	 In this work, our main goal is to prove that the solutions constructed in \cite{GMM17} converge strongly to solutions of the Navier-Stokes-Fourier system, at least when the initial data is in some of the spaces considered by the \cite{GMM17}. To do so, we draw inspiration from \cite{BMM} and decompose the Boltzmann equation into a system of two coupled equations, describing respectively the evolution of the macroscopic and microscopic parts of the initial data $f_\ini = \Pi f_\ini + f_{\ini, \bot}$.
	 \begin{itemize}
	 	\item One equation is posed in a ‘‘small Gaussian space" and its solution, generated by the macroscopic part $\Pi f_\ini$, is shown to converge strongly to a solution of the incompressible Navier-Stokes-Fourier system, by adapting the approach from \cite{BU91, IGT20}.
	 	\item The other equation is posed in a ‘‘large polynomial space" and describes the evolution of an initial layer generated by the microscopic part $f_{\ini, \bot}$.
	 \end{itemize}
	 We also prove that when the initial data $f_{\ini}$ is large, the corresponding solution to the Boltzmann equation exists at least up to the first singular time of the solution to the incompressible Navier-Stokes-Fourier system with initial data $f_{\ini, \WP}$.
	\subsection{Statement of the main result}
	
	Before we state the main result, let us define the functional spaces used in this paper.
	
	\begin{notation}
		\label{not:func_spc_main}
		For any positive Borelian function $m$, we denote by $L^p(m)$ the weighted Lebesgue space defined by the following norm:
		\begin{equation}
			\label{eqn:def_lebesgue_weight}
			\|f \|_{L^p\left( m\right)} := \|f m \|_{L^p},
		\end{equation} 
		and recall that~$M(v) = (2 \pi)^{-d/2} \exp(-|v|^2/2)$. We fix in this paper some $s > d/2$. For any $p \in [1, \infty]$ and $\alpha > 0$, we denote the norms
		\begin{gather*}
			\|f\|_\Spcp{p}{\alpha}:= \left\| f \langle v \rangle^\alpha\right\|_{L^p_v H^s_x},\\
			\|f\|_\Spcg{\beta} := \left\|f M^{-1/2} \langle v \rangle^\beta \right\|_{L^\infty_v H^s_x},
		\end{gather*}
		and corresponding functional spaces
		\begin{gather*}
			\Spcp{p}{\alpha} := \left\{ f=f(x, v) \,: \, \|f\|_\Spcp{p}{\alpha} < \infty \right\}, ~ p < \infty,\\
			\Spcp{\infty}{\alpha} := \left\{ f=f(x, v) \,: \, \|f\|_\Spcp{\infty}{\alpha} < \infty, \,  \|f(\cdot, v)\|_{H^s_x} \langle v \rangle^\alpha \xrightarrow[|v| \to \infty]{}  0 \right\},\\
			\Spcg{\beta} := \left\{ f=f(x, v) \,: \, \|f\|_\Spcg{\beta} < \infty, ~ \|f(\cdot, v)\|_{H^s_x} \langle v \rangle^\beta M^{-1/2} \xrightarrow[|v| \to \infty]{}  0 \right\}.
		\end{gather*}
		We also define the auxiliary norm
		\begin{equation*}
			\|h\|_{\Spcpp{p}{\alpha}}:= \left\| h \nu^{1/{p}}\right\|_{\Spcp{p}{\alpha}},
		\end{equation*}
		with $\nu$ being defined in \eqref{eqn:def_coll_freq}.
		Note that the assumption that the function vanishes for large $v$ in the case $p=\infty$ is made so that for any $\alpha > 0$ and $p_1 \geq 1$, the space~$\displaystyle \bigcap_{p \geq p_1} \Spcp{p}{\alpha}$ is dense in $\Spcp{\infty}{\alpha}$. Finally, for any real interval $I$ and Banach space~$X$, we denote $\CC_b(I ; X)$ the set of bounded continuous functions from $I$ to $X$.
	\end{notation}

	Our main result is the following.
	\begin{theo}
		\label{thm:main}
		Let $p \in [1, \infty]$ and $\beta > d/2 + 1$ (and recall that $s > d/2$), there exists some $\alpha_*(p) > 0$ such that, for any~$\alpha > \alpha_*(p)$, the following holds. Consider any initial data $f_\ini \in \Spcp{p}{\alpha}$, decomposed according to \eqref{eqn:initial_data_split}:
		\begin{gather*}
			f_\ini = f_{\ini, \WP} + f_{\ini, \IP} + f_{\ini, \bot},\\
			f_{\ini, \WP} (x, v) = \left( \rho_\ini^0(x) + u_\ini^0 (x) \cdot v+ \frac{1}{2} (|v|^2 - d) \theta_\ini^0(x)\right) M(v).
		\end{gather*}
		We make the following additional assumptions depending on the spatial domain:
		\begin{itemize}
			\item In the case $\Omega = \T^d$, we assume $(\rho_\ini^0, u^0_\ini, \theta_\ini^0)$ to be well-prepared, that is to say $\rho_\ini^0 + \theta_\ini^0 = 0$ and $\nabla_x \cdot u_\ini^0 = 0$. With the notations of Section \ref{scn:micro_macro_decomposition}, this means $\Pi f_\ini = f_{\ini, \WP}$, or equivalently $f_{\ini, \IP} = 0$. We also assume them to be mean-free.
			\item In the case $\Omega = \R^2$, we assume $(\rho_\ini^0, u^0_\ini, \theta_\ini^0) \in L^1_x$.
		\end{itemize}
		 Let~$\left(\rho^0, u^0, \theta^0\right)$ be the unique solution to \eqref{eqn:INSF} on a time interval $[0, T)$ (in the sense of Theorem \ref{thm:existence_INSF}) with initial data $(\rho^0_\ini, u^0_\ini, \theta^0_\ini)$, and define its kinetic counterpart
		\begin{equation*}
			f^0(t, x, v) := \left( \rho^0(t, x) + v \cdot u^0 (t, x) + \frac{1}{2} (|v|^2 - d) \theta_\ini^0(t, x)\right) M(v).
		\end{equation*}
		 There exists some positive $\eps_0=\eps_0(f_\ini, T)$ such that
		\begin{itemize}
			\item for any $\eps \in (0, \eps_0)$, the equation \eqref{eqn:bolt_scl} has a unique solution
			\begin{equation}
				\label{eqn:integ_solution}
				f^\eps \in \CC_b\left([0, T) ; \Spcp{p}{\alpha}\right) \cap L^p\left([0, T) ; \Spcpp{p}{\alpha}\right),
			\end{equation}
			\item the solution splits as
			\begin{gather*}
				f^\eps = f^0+ u^\eps_1 + u^\eps_\infty + u^\eps_\ac,\\
				u^\eps_1(0) = f_{\ini, \bot}, ~u^\eps_\infty(0) = 0, ~u^\eps_\ac(0) = f_{\ini, \IP},
			\end{gather*}
			where $u^\eps_1, u^\eps_\infty, u^\eps_\ac$ satisfy for some  $C, \gamma > 0$ depending only on $p, \alpha$ and $\beta$
			\begin{gather}
				\label{eqn:bound_u_1}
				\left\|u^\eps_1(t)\right\|_{\Spcp{p}{\alpha}} \leq C e^{-\gamma t / \eps^2} \|f_{\ini, \bot}\|_{\Spcp{p}{\alpha}},\\
				\label{eqn:vanish_u_infty}
				\lim_{\eps \to 0} \left(\sup_{0 \leq t < T} \left\|u^\eps_\infty(t)\right\|_{\Spcg{\beta}}\right) = 0,\\
				\label{eqn:vanish_u_w}
				u^\eps_\ac \rightharpoonup^* 0, \text{ in } L^\infty_t \Spcg{\beta},
			\end{gather}
			furthermore, if $\Omega = \R^d$, then $\left\|u^\eps_\ac\right\|_{L^q_t W_x^{s, \infty} \left( M^{-1/2} \langle v \rangle^\beta \right)} \to 0$ for any $q > \frac{2}{d-1}$.
			\item $u^\eps_\ac = 0$ if $f_{\ini, \IP} = 0$.
		\end{itemize}
	\end{theo}

	\begin{rem}
		In the physically relevant case $p = 1$, we have $\alpha_*(1) = 3$.
	\end{rem}

	\begin{rem}
		According to \cite[Theorem 1]{LZ01}, if the solution to \eqref{eqn:bolt_scl} is bounded in~$L^\infty_x L^1_v\left(\langle v \rangle^3\right)$, the solution to \eqref{eqn:boltz} is non-negative if $F_\ini$ is. This is indeed the case since $\Spcp{p}{\alpha} \subset \Spcp{1}{3}$ whenever $\alpha > \alpha_*(p)$, and $\Spcp{1}{3} \subset L^1_v L^\infty_x \left( \langle v \rangle^3 \right) \subset L^\infty_x L^1_v \left( \langle v \rangle^3 \right)$ because $s > d/2$.
	\end{rem}

	\begin{rem}
		Let us make some comments on the assumptions made on the initial data
		\begin{itemize}
			\item The mean-free assumption in the case of the torus is a compatibility condition coming from the fact that if $F^\eps = M + \eps f^\eps$ relaxes to the equilibrium $M$, then
			\begin{equation*}
				\int_{\R^d \times \T^d } f^\eps(t, x, v)
				\left(\begin{matrix}
					1\\v\\|v|^2
				\end{matrix}\right)
				\d v \d x \xrightarrow[t \rightarrow \infty]{}
				\left(\begin{matrix}
					0\\0\\0
				\end{matrix}\right),
			\end{equation*}
			and we recall that the following conservation laws hold:
			\begin{equation*}
				\int_{\R^d \times \T^d } f^\eps(t, x, v)
				\left(\begin{matrix}
					1\\v\\|v|^2
				\end{matrix}\right)
				\d v \d x = 
				\int_{\R^d \times \T^d } f_\ini( x, v)
				\left(\begin{matrix}
					1\\v\\|v|^2
				\end{matrix}\right)
				\d v \d x.
			\end{equation*}
			Thus, we need to assume
			\begin{equation*}
				\int_{\R^d \times \T^d } f_\ini(x, v)
				\left(\begin{matrix}
					1\\v\\|v|^2
				\end{matrix}\right)
				\d v \d x = 
				\left(\begin{matrix}
					0\\0\\0
				\end{matrix}\right).
			\end{equation*}
			\item The well-preparedness assumption is made because in the case $\Omega=\T^d$, the acoustic waves generated by the ill-prepared part do not disperse.
		\end{itemize}
	\end{rem}

	\begin{rem}
		Suppose the initial data is smooth (i.e. $H^{s+1}_x$). The part $u^\eps_\infty$ vanishes at a rate $\eps^{1/2}$, and if the macroscopic part of $f_\ini$ satisfies the incompressibility and Boussinesq conditions (i.e. $f_{\ini, \IP} = 0$ and thus $f_\ini = f_{\ini, \WP} + f_{\ini, \bot}$) then $u^\eps_\infty$ vanishes at a rate $\eps$. Also, the accoustic part $u^\eps_\ac$ is controled by $\eps t^{\frac{2}{d-1}}$ uniformly in $(x, v)$.
	\end{rem}

	\subsection{Notations}
	
	We will use the following notations throughout the rest of this paper. For any Banach spaces $X$ and $Y$, we will denote $\BBB(X;Y)$ the set of bounded linear operators from $X$ to $Y$, and abbreviate $\BBB(X)$ when $X=Y$. The relation denoted~$A \lesssim B$ is to be understood as $A \leq C B$ for some uniform constant $C > 0$. Finally, we will denote the convolution $f * g$ of maps $f : [0, \infty) \rightarrow X$, $g : [0, \infty) \rightarrow Y$ by
	$$f * g(t) = \int_0^t f(t-t') g(t') \d t'$$
	whenever the product $f(t_1) g(t_2)$ makes sense. For instance, when $f(t_1)$ is an operator acting on the vector $g(t_2)$, or when both are operators and we consider their composition. We introduce this definition so that, when the semigroup $S_L(t) = e^{t L}$ generated by a linear operator $L$ exists, the Duhamel principle applied to the evolution equation
	$$\partial_t u(t) = L u(t) + v(t).$$
	takes the simple form
	$$u(t) = S_L(t) u(0) + S_L*v(t)$$
	In particular, when all terms make sense, the following factorization formula holds for any linear operators $L, P$:
	\begin{equation*}
		S_{L+P} = S_L + S_L * \left(P S_{L+P}\right).
	\end{equation*}
	Furthermore, considering $L = (L + P) - P$, one also has $$S_L = S_{L+P} - S_{L+P} * \left(P S_{L}\right),$$
	and thus there holds
	\begin{equation}
		\label{eqn:semigroup_fact}
		S_{L+P} = S_L + S_L * \left(P S_{L+P}\right) = S_L + S_{L+P} * \left(P S_{L}\right).
	\end{equation}

	\subsection{Main reductions}
	\label{scn:main_reductions}
	
	We use an idea from \cite{BMM} and take advantage of the splitting of the linearized operator $\LL^\eps = \BB^\eps + \eps^{-2} \AA$ introduced in \cite{GMM17} and recalled in Section~\ref{scn:splitting_op}, where $\BB^\eps + \lambda / \eps^2$ generates a $\CC^0$-semigroup uniformly bounded in $\eps$, and there holds~$\AA \in \BBB\left( \Spcp{p}{\alpha} ; \Spcg{\beta} \right)$, for some $p \in [1, \infty]$ and positive~$\alpha, \lambda, \beta$. This allows to decompose the unknown of \eqref{eqn:bolt_scl} as $f^\eps(t) = h^\eps(t) + e^\eps(t) \in \Spcp{p}{\alpha} + \Spcg{\beta}$ where the parts $h^\eps$ and $e^\eps$ satisfy the following system of coupled equations:
	\begin{equation}
		\label{eqn:pre_coupled}
		\begin{cases}
			\partial_t h^\eps = \displaystyle \BB^\eps h^\eps + \frac{1}{\eps} Q\left( h^\eps, h^\eps \right) + 2 \eps^{-1} Q \left(h^\eps, e^\eps\right), ~ h^\eps_{| t = 0} = f_{\ini, \bot},\\
			\partial_t e^\eps = \displaystyle \LL^\eps e^\eps + \eps^{-1} Q\left( e^\eps, e^\eps \right) + \frac{1}{\eps^2} \AA h^\eps, ~ e^\eps_{| t = 0} = \Pi f_{\ini}.
		\end{cases}
	\end{equation}
	We will rewrite this system as a fixed point equation
	\begin{equation*}
		\Xi(h^\eps, e^\eps) = (h^\eps, e^\eps)
	\end{equation*}
	which will be solved using Banach's contraction theorem. We thus need to define properly $\Xi$ on some product space $\XX \times X$ and show contraction estimates for both coordinates $h^\eps$ and $e^\eps$. This means that we have to prove a priori estimates on both equations, which will dictate our choice for $\XX$ and $X$.
	
	The equation on $h^\eps$ will be studied using an energy method (Lemma \ref{lem:a_priori_poly}), which requires coercivity estimates on $\BB^\eps$ (Lemma \ref{lem:diss_B_eps}) and bounds for $Q$ (Lemma~\ref{lem:estimate_Q}).
	
	The equation on $e^\eps$ is — up to the coupling term $\eps^{-2} \AA h^\eps$ — the same equation as in \cite{BU91, IGT20} posed in the same functional space. We thus use the same approach and consider this equation in integral form:
	\begin{gather*}
		e^\eps = U^\eps \Pi f_{\ini} + \Psi^\eps\left( e^\eps, e^\eps \right) + \frac{1}{\eps^2} U^\eps * \AA h^\eps,
	\end{gather*}
	where $\Psi^\eps$ was defined in \eqref{eqn:def_psi_eps}, so as to rely on previously established results. Drawing inspiration from \cite{IGT20}, the part $e^\eps$ will be constructed indirectly: we need to identify a sub-part~$e_1^\eps$ such that there holds~$e^\eps = e_1^\eps + g^\eps$, where $g^\eps \to 0$ in $X$ as $\eps$ goes to $0$ and is what will actually be constructed during the fixed point argument.
	\begin{itemize}
		\item If one neglects the coupling term $\eps^{-2} U^\eps * \AA h^\eps$, previous works such as \cite{BU91, IGT20} suggest that $e^\eps$ should behave like $f^0 + U^\eps_\disp f_{\ini}$,
		where $f^0$ is the unique smooth solution to \eqref{eqn:KINSF} generated by $f_{\ini, \WP}$ on some time interval $[0, T)$ (in the sense of Theorem \ref{thm:existence_INSF}), and $U^\eps_\disp f_{\ini, \IP} = U^\eps_\disp f_\ini$ corresponds to acoustic waves oscillating with a rate of order $1/\eps$ (see Section \ref{scn:spectral_splitting} for the definition and properties of $U^\eps_\disp$).
		\item However, the coupling term $\eps^{-2} U^\eps * \AA h^\eps$ is not expected to vanish in $X$, which will be (a weighted in time version of) $L^\infty_t \Spcg{\beta}$, but satisfies (see Notations~\ref{not:convolution_split} for the definition of $\TT^\eps_1, \TT^\eps_\infty$ and Lemma \ref{lem:est_convo_1}-\ref{lem:est_convo_infty} for their respective estimates)
		\begin{gather*}
			\frac{1}{\eps^2} U^\eps * \AA h^\eps = \TT^\eps_1 h^\eps +  \TT^\eps_\infty h^\eps,\\
			\left\|\TT^\eps_q h^\eps \right\|_{L^q_t \Spcg{\beta}} \rightarrow 0, ~ (\eps \to 0), ~ q = 1, \infty.
		\end{gather*}
	\end{itemize}
	We therefore define the sub-part mentioned above as $e_1^\eps = f^0 + U^\eps_\disp f_{\ini} + \TT^\eps_1 h^\eps$, thus the pair $(h^\eps, g^\eps)$ solves the system
	\begin{equation}
		\label{eqn:coupled_diff}
		\begin{cases}
			\partial_t h^\eps &= \displaystyle \BB^\eps h^\eps + \frac{1}{\eps} Q\left( h^\eps, h^\eps \right) + \frac{2}{\eps} Q \left(h^\eps, g^\eps + f^0 + \TT^\eps_1 h^\eps + U^\eps_{\disp} f_{\ini}\right),\\
			h^\eps(0) &= f_{\ini, \bot},\\
			g^\eps &= \SS^\eps[h^\eps] + \Psi^\eps(g^\eps, g^\eps) + \Phi^\eps[h^\eps] g^\eps,
		\end{cases}
	\end{equation}
	where we have denoted
	\begin{gather*}
		\SS^\eps[h^\eps] := \SS^\eps_0 + \Psi^\eps\left( \TT^\eps_1 h^\eps + U^\eps_\disp f_{\ini} , 2 f^0 + \TT^\eps_1 h^\eps + U^\eps_\disp f_{\ini} \right) + \TT^\eps_\infty h^\eps,\\
		\SS^\eps_0 := \left(U^\eps - U^\eps_\disp - U^0\right) \Pi f_{\ini} + \left(\Psi^\eps - \Psi^0\right)\left(f^0, f^0\right),\\
		\Phi^\eps[h^\eps] := 2 \Psi^\eps\left(  f^0 + \TT^\eps_1[h^\eps] + U^\eps_\disp f_{\ini}, \, \cdot \, \right),
	\end{gather*}
	and we used the fact that $U^\eps_\disp = U^\eps_\disp \Pi$ (see Section \ref{scn:spectral_splitting}).
	Note that $f_\ini$, $f^0$ and $T$ are fixed so we do not indicate them in the notations above.
	
	\subsection{Comparison with previous works}
	
	In \cite{BMM}, the authors considered the system~\eqref{eqn:pre_coupled} in order to derive estimates on both parts uniformly in $\eps$. They choose to split the initial data as $h^\eps(0) = f_\ini$ and $e^\eps(0) = 0$, which allowed to rely on hypocoercivity results from \cite{Briant} valid for $e^\eps(0)$ small enough. Like in the present work, they obtained the following control for the coupling term: $\eps^{-2} \AA h^\eps(t) = O\left( \eps^{-2} e^{-\lambda t / \eps^2} \right)$ which leads to a term of order $O(1)$ in Gronwall estimates.
	
	Here, we need to show that $e^\eps$ is not only bounded but converges strongly to~$f^0$ as well when~$\eps$ goes to zero. To this end, we use the integral formulation of the equation on $e^\eps$ which allows to use the spectral study of the linearized operator~$\LL^\eps$ led in such works as \cite{EP75, U86, G21, YY16}, this is the approach used in \cite{BU91, IGT20}.
	
	Note that if $f_\ini$ were purely macroscopic, then it would belong to $\Spcg{\beta}$ and \eqref{eqn:pre_coupled} would reduce to the equation considered in \cite{BU91, IGT20}. The latter differs from the equation on $e^\eps$ in \eqref{eqn:pre_coupled} by the coupling term $\eps^{-2} U^\eps * \AA h^\eps$, which may not be small in $L^\infty_t \Spcg{\beta}$, but is small in $L^\infty_t \Spcg{\beta} + L^1_t \Spcg{\beta}$. This is made possible (1) by the splitting of the initial data $h^\eps(0) = f_{\ini, \bot}$, (2) by generalizing some properties of $U^\eps$ known to hold in $\Spcg{\beta}$ to the larger space $\Spcp{p}{\alpha}$, using the theory of space enlargement from \cite{M05, MM16, GMM17}.
	
	\subsection{Outline of the paper}
	
	In Section \ref{scn:study_polynomial}, we prove coercivity estimates for $\BB^\eps$ and bounds for $Q$ which we use to prove the well-posedness of the following equation for given source terms $h$ and $g$
	\begin{equation*}
		\partial_t f = \BB^\eps f + \frac{1}{\eps} Q(h, h + g).
	\end{equation*}
	In particular, we prove the stability estimate required for the mapping $\Xi$ from Section~\ref{scn:main_reductions} to be a contraction.

	In Section \ref{scn:asymptotics_semigroup}, we recall the spectral properties of $\LL + \eps v \cdot \nabla_x$ from \cite{EP75, U86, G21, IGT20}, which dictate the asymptotics of the semigroup $U^\eps$ in some Gaussian space. We then extend some asymptotic properties to polynomial spaces $\Spcp{p}{\alpha}$ using space enlargement theory.
	
	In Section \ref{scn:a_priori_gauss}, we prove the necessary estimates to solve
	\begin{equation*}
		g = \Psi^\eps(g, g) + \Phi^\eps[h] g + \SS^\eps[h]
	\end{equation*}
	for some unique small $g$ when $\eps$ is small enough using a contraction argument. To be more precise, we show that $\Psi^\eps$ is uniformly bounded in $\eps$, that $\Phi^\eps[h]$ is a contraction with Lipschitz constant depending on a bound of $h$ and for an equivalent norm depending on~$f^0$, and that $\SS^\eps[h]$ vanishes as $\eps$ goes to zero.
	
	\section{Study in the polynomial space}
	
	\resetcounter
	\label{scn:study_polynomial}
	
	\subsection{Splitting of the linearized operator}
	\label{scn:splitting_op}
	
	In \cite{GMM17}, the authors present a splitting of the linearized Boltzmann operator $\LL$, defined as \eqref{eqn:def_L}, in the form
	\begin{equation*}
		\LL = \AA_\delta + \BB_\delta = -\nu + \overline{\BB}_\delta + \AA_\delta, ~ \delta \in (0, 1),
	\end{equation*}
	and $\AA_\delta$, $\overline{\BB}_\delta$ are defined by
	\begin{gather*}
		\AA_\delta f(v) := \int_{\R^d \times \mathbb{S}^{d-1} } \Theta_\delta \left(M'_* f' + M' f'_* - M f_*\right) |v-v_*| \d v_* \d\sigma,\\
		\overline{\BB}_\delta f(v)  := \int_{\R^d \times \mathbb{S}^{d-1} } \left(1-\Theta_\delta\right) \left(M'_* f' + M' f'_* - M f_*\right) |v-v_*| \d v_* \d\sigma,
	\end{gather*}
	where we recall that $\nu$ is defined by \eqref{eqn:def_coll_freq}, the notations $v', v'_*, g', g'_*, g_*$ are defined in Section \ref{scn:boltzmann}, the function $\Theta_\delta=\Theta_\delta(v, v_*, \sigma)$ is smooth, bounded by one on
	$$\{ |v| \leq \delta^{-1}, 2 \delta \leq |v-v_*| \leq \delta^{-1}, |\cos \theta | \leq 1-2 \delta\},$$
	and supported in
	$$\{ |v| \leq 2\delta^{-1}, \delta \leq |v-v_*| \leq 2\delta^{-1}, |\cos \theta | \leq 1- \delta\},$$
	where $\cos \theta := \sigma \cdot (v-v_*)/|v-v_*|$.
		
	\medskip
	For any $p \in [1, \infty]$, $\alpha \geq 0$, $\beta \geq 0$, and $\delta \in (0, 1)$, the operator $\AA_\delta$ is bounded from~$L^p\left( \langle v \rangle^\alpha \right)$ to $L^\infty\left( \langle v \rangle^\beta M^{-1/2} \right)$ as can be seen from Carleman's representation:
	\begin{equation*}
		\AA_\delta f(v) = \int_{\R^d} k_\delta (v, v_*)f(v_*) \d v_*,
	\end{equation*}
	where $k_\delta \in \CC^\infty_c(\R^d \times \R^d)$ (see \cite[(4.9)]{GMM17}). By its integral nature in $v$, it is clear that
	\begin{equation}
		\label{eqn:reg_AA}
		\AA_\delta : \Spcp{p}{\alpha} \rightarrow \Spcg{\beta}
	\end{equation}
	is also bounded. In the rest of the paper, we will be interested in the scaled version of $\BB_\delta$ together with a transport term, which we denote
	\begin{equation}
		\label{eqn:def_B_eps}
		\BB^\eps := \frac{1}{\eps^2}\left(\BB_\delta + \eps v \cdot \nabla_x\right).
	\end{equation}
	
	\begin{lem}
		\label{lem:diss_B_eps}
		Let $p \in [1, \infty)$, there exists $\alpha_\BB(p) > 2$ such that for any $\alpha > \alpha_\BB(p)$, if $\delta$ is small enough, $\overline{\BB}_\delta$ is $\nu$-bounded with relative bound less than one, and there holds for any $\eps \in (0, 1)$ and $h \in \Spcp{p}{\alpha}$
		\begin{gather}
			\label{eqn:coercivity}
			\int_{\R^d} \langle \BB^\eps h, h \rangle_{H^s_x} \|h\|_{H^s_x}^{p-2} \langle v \rangle^{p \alpha}\d v \leq -\frac{\sigma_\BB(p, \alpha)}{\eps^2} \|h\|^p_{\Spcpp{p}{\alpha}},\\
			\notag
			\sigma_\BB(p, \alpha)=\sigma_\BB(p, \alpha ; \delta) > 0,
		\end{gather}
		where $\langle \cdot, \cdot \rangle_{H^s_x}$ denotes the usual inner product in $H^s_x$, and $\Spcpp{p}{\alpha}$ was defined in Notations \ref{not:func_spc_main}.
	\end{lem}

	\begin{proof}
		In \cite{GMM17}, the authors prove a control on~$ \overline{\BB}_\delta$ in the space $\Spcp{p}{\alpha}$ by considering its positive variant:
		\begin{gather*}
			\widetilde{\BB}_\delta f(v)  := \int_{\R^d \times \mathbb{S}^{d-1} } \left(1-\Theta_\delta\right) \left(M'_* f' + M' f'_* + M f_*\right) |v-v_*| \d v_* \d\sigma,
		\end{gather*}
		so that one has by the triangle inequality
		\begin{equation}
			\label{eqn:triangle_b_delta}
			\left\| (\overline{\BB}_\delta h) (v)\right\|_{H^s_x} \leq \bigg(\widetilde{\BB}_\delta \left(\|h\|_{H^s_x}\right) \bigg)(v),
		\end{equation}
	 	and by proving the following estimate\footnote{This estimate is actually proved for $\overline{\BB}_\delta$ in \cite[(4.32)]{GMM17} , but as the authors point out in page~74, the same estimate holds for $\widetilde{\BB}_\delta$ with the same proof.} for any $\alpha > 2$, where $o_\delta(1)$ designates any quantity that vanishes as $\delta \to 0$:
	 	\begin{gather}
	 		\label{eqn:bound_b_delta_tilde}
	 		\left\| \widetilde{\BB}_\delta f \right\|_{L^p\left( \langle v \rangle^\alpha \right)} \leq \left(\phi_p(\alpha) + o_{\delta}(1)\right) \|\nu f\|_{ L^p\left( \langle v \rangle^\alpha \right) },\\
	 		\notag
	 		\phi_p(\alpha) := \frac{4}{(\alpha + 2)^{1/p}(\alpha - 1)^{1/{p'}}}.
	 	\end{gather}
 		Let us now prove Lemma \ref{lem:diss_B_eps}: by the definition of $\BB_\delta$ and the skew-adjointness of~$v \cdot \nabla_x$ in $H^s_x$ for any $v \in \R^d$, one has
 		\begin{align*}
			\int_{\R^d} \big\langle \BB_\delta h+ &\eps v \cdot \nabla_x h, h \big\rangle_{H^s_x} \|h\|_{H^s_x}^{p-2} \langle v \rangle^{p \alpha}\d v\\
			& = - \|h\|_{\Spcpp{p}{\alpha}}^p + \int_{\R^d} \big\langle \overline{\BB}_\delta h, h \big\rangle_{H^s_x} \|h\|_{H^s_x}^{p-2} \langle v \rangle^{p \alpha}\d v \\
			& \leq - \|h\|_{\Spcpp{p}{\alpha}}^p + \int_{\R^d} \|\overline{\BB}_\delta h\|_{H^s_x} \|h\|_{H^s_x}^{p-1} \langle v \rangle^{p \alpha}\d v.
 		\end{align*}
 		Using \eqref{eqn:triangle_b_delta} and denoting to lighten the notation $g(v) := \|\nu^{1/p} h(v)\|_{H^s_x}$, we have
 		\begin{align*}
 			\int_{\R^d} \|\BB_\delta h\|_{H^s_x} &\|h\|_{H^s_x}^{p-1} \langle v \rangle^{p \alpha}\d v \\
 			&  \leq \int_{\R^d} \left(\nu^{-1/{p'}} \widetilde{\BB}_\delta \nu^{-1/p} \right) (g) \times g^{p-1} \langle v \rangle^{p \alpha}\d v \\
 			& \leq \left\| \nu^{-1/{p'}} \widetilde{\BB}_\delta \nu^{-1/p}  \right\|_{ \BBB(L^p\left( \langle v \rangle^\alpha \right)) } \|g\|_{L^p\left( \langle v \rangle^\alpha \right)}^{p}.
 		\end{align*}
 		Using the estimate $\nu_0 \langle v \rangle \leq \nu(v) \leq \nu_1 \langle v \rangle$ and \eqref{eqn:bound_b_delta_tilde}, one deduces that, under the condition $\alpha > 2 + 1/{p'}$, there holds
 		\begin{align*}
 			\left\| \nu^{-1/{p'}} \widetilde{\BB}_\delta \nu^{-1/p}  \right\|_{ \BBB(L^p\left( \langle v \rangle^\alpha \right)) } =& \left\| \nu^{-1/{p'}} \left(\widetilde{\BB}_\delta \nu^{-1}\right) \nu^{1/q}\right\|_{\BBB(L^p\left( \langle v \rangle^\alpha \right)) } \\
 			\leq& \left\| \nu^{-1/{p'}} \right\|_{ \BBB\left(L^p\left(\langle v \rangle^{\alpha - 1/{p'}}\right) ;  L^p\left(\langle v \rangle^{\alpha}\right)\right) } \\
 			&\times \left\| \widetilde{\BB}_\delta \nu^{-1} \right\|_{ \BBB\left(L^p\left(\langle v \rangle^{\alpha - 1/{p'}}\right) \right) }\\
 			& \times \left\| \nu^{1/p} \right\|_{ \BBB\left(L^p\left(\langle v \rangle^{\alpha}\right) ; L^p\left(\langle v \rangle^{\alpha - 1/{p'}}\right) \right) }\\
 			\leq& \left(\frac{\nu_1}{\nu_0}\right)^{1/{p'}} \bigg( \phi_p(\alpha + 1/{p'}) + o_\delta(1) \bigg).
 		\end{align*}
 		To sum up, we have shown 
 		\begin{align*}
 			\int_{\R^d} \langle \BB^\eps h, h \rangle_{H^s_x} &\|h\|_{H^s_x}^{p-2} \langle v \rangle^{p \alpha}\d v \\
 			& \leq -\frac{1}{\eps^2}\bigg(1 -  \left(\frac{\nu_1}{\nu_0}\right)^{1/{p'}} \phi_p(\alpha + 1 /{p'}) + o_\delta(1)\bigg)\|h\|_{\Spcpp{p}{\alpha}}^p\\
 			& =: - \frac{\sigma_\BB(p, \alpha)}{\eps^2} \|h\|_{\Spcpp{p}{\alpha}}^p.
 		\end{align*}
 		The quantity $\sigma_\BB(p, \alpha)$ is positive as soon as $\alpha > \alpha_\BB(p)$ and $\delta$ is small enough, where we defined
		\begin{gather}
 			\label{eqn:def_alpha_B}
 			\alpha_\BB(p) := \inf \left\{ \alpha > 2 + 1/p' ~ : ~ \left(\frac{\nu_1}{\nu_0}\right)^{1/{p'}} \phi_p(\alpha + 1/{p'}) < 1\right\},
 		\end{gather}
 		which concludes the proof of the lemma.
	\end{proof}

	\begin{rem}
		Note that the definition of $\sigma_\BB(p, \alpha)$ actually depends on the choice of~$\delta$, but this does not make any difference in the rest of the paper as long as $\sigma_\BB(p, \alpha)=\sigma_\BB(p, \alpha ; \delta) > 0$, so we do not mention it from now on as to lighten the notations.
	\end{rem}

	\begin{rem}
		The threshold $\alpha_\BB(p)$ is monotonically increasing in $p$, and in the cases~$p = 1$ and $p = \infty$, this threshold is simply
		\begin{equation*}
			\alpha_\BB(p) =
			\begin{cases}
				2, &\text{ when } p = 1,\\
				\displaystyle 1 + 4 \, \frac{\nu_1}{\nu_0} > 5, &\text{ when } p = \infty.
			\end{cases}
		\end{equation*}
	\end{rem}

	\subsection{A priori estimates}

	In this section, we present the estimates necessary to the a priori study of the equation on $h^\eps$ in \eqref{eqn:coupled_diff}. We first prove a bound for the collision operator $Q$ (Lemma \ref{lem:estimate_Q}), define for functions $t \mapsto h(t) \in \Spcp{p}{\alpha}$ a norm which measures the exponential decay and moment gain induced by $\BB^\eps$  in \eqref{eqn:coupled_diff} uniformly in $\eps$ reminiscent of the ones from \cite[Section 5.1]{ALT}, and finally we give a priori estimates for the equation on $h^\eps$ in \eqref{eqn:coupled_diff} (Lemma \ref{lem:a_priori_poly}).
	
	\medskip
	This control on $Q$ is inspired of \cite[Lemma 2.3]{HTT} in which the case $p = 2$ is treated for hard potentials without cut-off.
	
	\newcommand{\Integ}{\mathbf{I}}
	\begin{lem}
		\label{lem:estimate_Q}
		Let $p \in [1, \infty]$, there exists $\alpha_Q(p) > 2$ such that, for any $\alpha > \alpha_Q(p)$, there exists some $C=C(p, \alpha) > 0$ satifying
		\begin{equation}
			\label{eqn:bound_Q}
			\left\| Q(f, g) \right\|_{\Spcp{p}{\alpha - 1/{p'} }} \leq C \left(\|f\|_{\Spcp{p}{\alpha}} \|g\|_{\Spcpp{p}{\alpha}} + \|f\|_{\Spcpp{p}{\alpha}} \|g\|_{\Spcp{p}{\alpha}}\right),
		\end{equation}
		and if $p < \infty$,
		\begin{align}
			\label{eqn:dual_bound_Q}
			\int_{\R^d} \left\langle Q(f, g), h\right\rangle_{H^s_x} &\|h\|_{H^s_x}^{p-2} \langle v \rangle^{p \alpha} \d v\\
			\notag
			& \leq C \left(\|f\|_{\Spcp{p}{\alpha}} \|g\|_{\Spcpp{p}{\alpha}} + \|f\|_{\Spcpp{p}{\alpha}} \|g\|_{\Spcp{p}{\alpha}}\right) \|h\|_{\Spcpp{p}{\alpha}}^{p-1},
		\end{align}
	\end{lem}

	\begin{proof}
		In this proof, we denote $q := p'$ so as not to cause confusion with the pre-collisional velocity notations $v', v'_*$, and we define the threshold
		\begin{gather}
			\label{eqn:def_alpha_Q}
			\alpha_Q(p) := 2 + \frac{d}{q}.
		\end{gather}
		\step{1}{Main reductions}
		First, note that the algebra structure of $H^s_x$ implies
		\begin{align*}
			\left|\langle Q^+(f(v), g(v)), h(v) \rangle_{H^s_x}\right| & = \left| \int_{\R^d} |v-v_*| \langle f(v') g(v'_*), h(v) \rangle_{H^s_x} \d v_* \right| \\
			& \leq \|h(v)\|_{H^s_x}  \int_{\R^d} |v-v_*| \|f(v') g(v'_*)\|_{H^s_x} \d v_*\\
			& \lesssim \|h(v)\|_{H^s_x} \, Q^+\left( \|f(v)\|_{H^s_x}, \|g(v)\|_{H^s_x} \right).
		\end{align*}
		The same estimate holds for $Q^-$, thus, denoting
		\begin{gather*}
			F(v) := \|f(v)\|_{H^s_x}, ~ G(v) := \|g(v)\|_{H^s_x},\\
			m(v) := \langle v \rangle^\alpha,
		\end{gather*}
 		it is enough to prove
		\begin{align*}
			\left| \int_{\R^d} Q(F, G) H \,  m^r \, \d v  \right| \leq C \bigg(&\|F\|_{L^p\left(\langle v \rangle^{1/p} m\right)} \|G\|_{L^p\left(m\right)} \\
			&+ \|F\|_{L^p\left(m\right)} \|G\|_{L^p\left(\langle v \rangle^{1/p} m\right)}\bigg) \left\| H m^{r-1} \langle v \rangle^{1/q}\right\|_{L^q},
		\end{align*}
		where we recall the notation $L^p(m)$ is defined in \eqref{eqn:def_lebesgue_weight}. Indeed, this control will imply
		\begin{itemize}
			\item \eqref{eqn:bound_Q} for~$H(v) := \|h(v)\|_{H^s_x}$ and $r=2$, by duality,
			\item \eqref{eqn:dual_bound_Q} for $H(v) := \|h(v)\|_{H^s_x}^{p-1}$ and $r = p$.
		\end{itemize}
		To do so, we only need to estimate the integral coming from $Q^-$:
		\begin{align*}
			\left|\int_{\R^d} Q^-(F, G) H m^r \d v\right| \lesssim \,& \int |v-v_*| |F_*| |G| |H| m^r \d v \d v_*
			=: \Integ_1,
		\end{align*}
		and the one coming from $Q^+$. For the latter, we use the following estimate from \cite[(2.1.15)]{AMUXY}:
		\begin{equation*}
			|m-m'| \lesssim \sin\left(\frac{\theta}{2}\right) \left(m' + \left\langle v'_*\right\rangle \left\langle v'\right\rangle^{\alpha-1} + \sin\left(\frac{\theta}{2}\right)^{\alpha-1} m'_* \right),
		\end{equation*}
		where we recall the notations $m' := m(v')$, $m_*:= m(v_*)$ and $m'_* := m(v_*')$. We will use it under the simpler form
		\begin{equation*}
			m \lesssim m' + \langle v_*' \rangle \langle v' \rangle^{\alpha - 1} + m_*' \theta^\alpha,
		\end{equation*}
		which suggests the splitting
		\begin{align*}
			\left|\int_{\R^d} Q^+(F, G) H m^r \d v\right| \lesssim\,& \int |v- v_*| |F'_*| |G'| |H| m' m^{r-1} \d v \d v_* \d\sigma \\
			& +   \int |v- v_*| |F'_*| |G'| |H| \langle v_*' \rangle \langle v' \rangle^{\alpha-1} m^{r-1} \d v \d v_* \d\sigma\\
			& +  \int |v- v_*| |F'_*| |G'| |H| \theta^\alpha m_*' m^{r-1} \d v \d v_* \d\sigma \\
			= :&\Integ_{2} +  \Integ_{3} + \Integ_{4}.
		\end{align*}

		\step{1}{Estimating $\Integ_1$ and $\Integ_3$}
		We start by using the fact that
		$$|v-v_*| \leq \langle v \rangle \langle v_* \rangle \leq \langle v \rangle^{1/p} \langle v \rangle^{1/q} \langle v_* \rangle,$$
		which implies the following for $\Integ_1$:
		\begin{align}
			\notag
			\Integ_1 & \lesssim \int \left(|G| m \langle v \rangle^{1/p}\right) \left(|H| m^{r-1} \langle v \rangle^{1/q}\right) \d v \int |F_*| \langle v_*\rangle \d v_*\\
			\label{eqn:est_tri}
			& \lesssim \| G \|_{L^p\left( m \langle v \rangle^{1/p} \right)} \left\| H m^{r-1} \langle v \rangle^{1/q}\right\|_{L^q} \|F\|_{L^1 \left( \langle v \rangle \right) }
		\end{align}
		where we used Hölder's inequality in the last line. We deal with~$\Integ_{3}$ in a similar fashion, using $|v-v_*|=|v'-v'_*|= \leq \langle v_*' \rangle \langle v' \rangle$, which yields this time
		$$\Integ_{3} \lesssim \| G \|_{L^p\left( m \langle v \rangle^{1/p} \right)} \left\| H m^{r-1} \right\|_{L^q} \|F\|_{L^1 \left( \langle v \rangle^2 \right) }.$$
		
		\step{2}{Estimating $\Integ_2$}
		We use the following identity, which is a direct consequence of \eqref{eqn:pre_collision_velocities}
		\begin{equation*}
			|v - v_*| = |v'-v_*'| = \sqrt{\frac{1 - \cos \theta}{2}} |v-v_*'|
		\end{equation*}
		and thus
		\begin{align*}
			|v-v_*| \leq |v'-v_*'|^{1/p} |v-v'_*|^{1/q} \leq \langle v' \rangle^{1/p} \langle v \rangle^{1/q} \langle v'_* \rangle.
		\end{align*}
		From this and Hölder's inequality, we get the following control on $\Integ_{2}$:
		\begin{align*}
			\Integ_{2} \leq& \left[\int \bigg(\langle v_*' \rangle^{1/p} |F_*'|^{1/p} |G'| \langle v' \rangle^{1/p} m'\bigg)^p \d v \d v_* \d \sigma\right]^{1/p}\\
			& \times \left[  \int \bigg( \langle v'_* \rangle^{1/q} |F_*'|^{1/q} |H| \langle v \rangle^{1/q} m^{r-1} \bigg)^q \d v \d v_* \d \sigma  \right]^{1/q}.
		\end{align*}
		Using the change of variables $(v, v_*) \mapsto (v', v'_*)$ in the first integral and $v_* \mapsto v'_*$ in the second, whose respective Jacobian determinants are
		$$\left| \frac{\d (v, v_*)}{\d (v', v'_*)} \right| = 1, \quad \left|\frac{\d v_*'}{\d v_*}\right| = \frac{1 + \cos \theta}{8},$$
		we obtain the estimate
		\begin{align*}
			\Integ_{2}  \leq &\left[\int \bigg(\langle v_* \rangle^{1/p} |F_*|^{1/p} |G| \langle v \rangle^{1/p} m\bigg)^p \d v \d v_*\right]^{1/p}\\
			& \times \left[  \int \bigg( \langle v_* \rangle^{1/q} |F_*|^{1/q} |H| \langle v \rangle^{1/q} m^{r-1} \bigg)^q \d v \d v_* \d \sigma \right]^{1/q}.
		\end{align*}
		We conclude that $\Integ_{2}$ satisfies \eqref{eqn:est_tri} using again Hölder's inequality on each integral.
		
		\step{3}{Estimating $\Integ_4$}
		We use this time the bound
		$$|v' - v| = \sqrt{\frac{1 - \cos \theta}{2}} |v-v_*| \lesssim \theta^{-1} |v-v_*|,$$
		which yields
		$$|v'-v| \lesssim \theta^{-1/q} |v-v'|^{1/q} |v'-v'_*|^{1/p} \lesssim \theta^{-1/q} \langle v' \rangle \langle v'_* \rangle^{1/p} \langle v \rangle^{1/q}.$$
		Plugging this estimate in $\Integ_{4}$, we get
		\begin{align*}
			\Integ_{4} \lesssim& \int \theta^{\alpha - (1 + 1/q)} \left( \langle v \rangle |g| \right)' \left( |F| m \langle v \rangle^{1/p} \right)'_* \left( |H| m^{r-1} \langle v \rangle^{1/q} \right) \d v \d v_* \d \sigma\\
			\lesssim &\Bigg(\int \theta^{\alpha - (1 + 1/q)} \left( \langle v \rangle |g| \right)' \bigg[\left( |F| m \langle v \rangle^{1/p} \right)'_* \bigg]^p \d v \d v_* \d \sigma \Bigg)^{1/p}\\
			& \times \Bigg(\int \theta^{\alpha - (1 + 1/q)} \left( \langle v \rangle |G| \right)' \left( |H| m^{r-1} \langle v \rangle^{1/q} \right)^q  \d v \d v_* \d \sigma \Bigg)^{1/q}
		\end{align*}
		where we used Hölder's inequality. Using then $(v, v_*) \mapsto (v', v'_*)$ in the first integral and $v' \mapsto v_*$ in the second one, which is such that $\displaystyle \left|\frac{\d v_*}{\d v'}\right| \lesssim \theta^{-2}$, one gets
		\begin{align*}
			\Integ_{4} &\lesssim \|G\|^{1/p}_{L^1\left( \langle v \rangle \right)} \|F\|_{L^p\left( m \langle v \rangle^{1/p} \right)} \times \Bigg(\int \theta^{\alpha - (3 + 1/q)} \left( \langle v \rangle |G| \right)' \left( |H| m^{r-1} \langle v \rangle^{1/q} \right)^q \d v \d v' \d \sigma\Bigg)^{1/q}.
		\end{align*}
		Recall that $\theta$ is the angle formed by $(v-v_*, \sigma)$, we have to replace it by the angle~$\psi$ formed by $(v-v', \sigma)$ before rewriting the previous integral factor in spherical coordinates. The following relations hold:
		\begin{gather*}
			\psi = \frac{\pi - \theta}{2}, \quad
			\cos \psi = \frac{v-v'}{|v-v'|} \cdot \sigma.
		\end{gather*}
		We then bound $\Integ_{4}$ as
		\begin{align*}
			 \Integ_{4} \lesssim \bigg( \int (\pi - 2 \psi)^{\alpha - (3+1/q) } \sin \psi^{d-2} \d \psi\bigg)  \times  \|G\|_{L^1\left( \langle v \rangle \right)} \|F\|_{L^p\left( m \langle v \rangle^{1/p} \right)} \left\| Hm^{r-1} \langle v \rangle^{1/q}\right\|_{L^q}
		\end{align*}
		Because of the assumption on $\alpha$, we have $\alpha - (3 + 1/q) > -1$, therefore the integral factor is finite. Finally, we have that, for some positive $C_0=C_0(\alpha)$,
		\begin{align*}
			\Integ_j &\leq C_0 \| G \|_{L^p\left( m \langle v \rangle^{1/p} \right)} \left\| H m^{r-1} \langle v \rangle^{1/q}\right\|_{L^q} \|F\|_{L^1 \left( \langle v \rangle^2 \right) }\\
			&\leq C_0 \|\langle v \rangle^{(2-\alpha)q}\|_{L^1}^{1/q} \| G \|_{L^p\left( m \langle v \rangle^{1/p} \right)} \left\| H m^{r-1} \langle v \rangle^{1/q}\right\|_{L^q} \|F\|_{L^p \left( m \right) }.
		\end{align*}
		The conclusion follows with $C=C_0 \|\langle v \rangle^{(2-\alpha)q}\|_{L^1}^{1/q}$, which is finite because the assumption on $\alpha$ implies $(2 - \alpha)q < -d$.
	\end{proof}
	
	We now prove the a priori estimates for the equation on $h^\eps$ from \eqref{eqn:coupled_diff}, but first let us introduce a norm which measures the exponential decay and the moment gain induced by $\BB^\eps$.
	
	\begin{notation}
		\label{not:func_spc}
		For any $p \in [1, \infty]$, $\alpha > \alpha_\BB$, $\sigma \in \big(0, \sigma_\BB\big)$ and $\eps \in (0, 1)$, where $\alpha_\BB$ and $\sigma_\BB$ are defined in Lemma \ref{lem:diss_B_eps}, denote the functional spaces
		\begin{equation*}
			\tSpcp{p}{\alpha} = \CC_b\left( [0, T) ; \Spcp{p}{\alpha} \right) \cap L^p\left( [0, T) ; \Spcpp{p}{\alpha} \right),
		\end{equation*}
		where $\CC_b$ is defined in Notations \ref{not:func_spc_main}, $T$ in Section \ref{scn:main_reductions}, and endowed with the norm
		\begin{gather*}
			\tSpcpN{h}{p}{\alpha}^p:= \sup_{0 \leq t \leq T} \Bigg( \frac{1}{p} \left(e^{\sigma t / \eps^2} \|h(t)\|_\Spcp{p}{\alpha}\right)^p + \frac{\sigma_\BB - \sigma}{\eps^2} \int_{0}^{t} \left(e^{\sigma t' / \eps^2} \|h(t')\|_\Spcpp{p}{\alpha}\right)^p \d t' \Bigg),
		\end{gather*}
		if $p < \infty$, and
		$$\tSpcpN{h}{\infty}{\alpha} := \sup_{0 \leq t \leq T} \left\|e^{\sigma t / \eps^2} h(t)\right\|_{\Spcp{\infty}{\alpha}}$$
		otherwise.
	\end{notation}
	
	\begin{lem}
		\label{lem:a_priori_poly}
		Let $\alpha > \max\{\alpha_\BB, \alpha_Q\}$, $\sigma \in (0, \sigma_\BB)$, $\beta \geq 0$.  For any $h \in L^\infty\left([0, T) ; \Spcp{p}{\alpha}\right) \cap L^p\left([0, T) ; \Spcpp{p}{\alpha}\right)$ and $g \in L^\infty([0, T) ; \Spcg{\beta})$, the evolution equation
		\begin{equation*}
			\begin{cases}
				\partial_t \solp = \displaystyle \BB^\eps \solp + \frac{1}{\eps} Q(h, h)  + \frac{1}{\eps} Q(h, g),\\
				\solp(0) \in \Spcp{p}{\alpha},
			\end{cases}
		\end{equation*}
		has a unique solution $\solp \in \tSpcp{p}{\alpha}$. Furthermore, it satisfies for some $C=C(p, \alpha, \beta)$
		\begin{equation}
			\label{eqn:a_priori_poly}
			\tSpcpN{\solp}{p}{\alpha} \leq  C \eps \tSpcpN{h}{p}{\alpha} \bigg(\tSpcpN{h}{p}{\alpha} + \|g\|_{L^\infty_t \Spcg{\beta}}\bigg) + \|\solp(0)\|_\Spcp{p}{\alpha},
		\end{equation}
		as well as the stability estimate for any pair of solutions $\overline{h}_1, \overline{h}_2$ and corresponding source terms $h_1, g_1, h_2, g_2$:
		\begin{align}
			\notag
			\tSpcpN{\solp_1-\solp_2}{p}{\alpha}\leq&  C \eps \tSpcpN{h_1 - h_2}{p}{\alpha} \bigg(\tSpcpN{h_1+h_2}{p}{\alpha} +  \|g_1\|_{L^\infty_t \Spcg{\beta}}\bigg) \\
			\label{eqn:a_priori_poly_stab}
			&+ C \eps \tSpcpN{h_2}{p}{\alpha} \|g_1 - g_2\|_{L^\infty_t \Spcg{\beta}}
			 + \|\solp_1(0) - \solp_2(0)\|_\Spcp{p}{\alpha}.
		\end{align}
	\end{lem}

	\begin{proof}
		We will denote in this proof $\tspcp = \tSpcp{p}{\alpha}$, $\spcp = \Spcp{p}{\alpha}$, $\spcpp = \Spcpp{p}{\alpha}$ to lighten the notations. As the constant $C$ in Lemma \ref{lem:estimate_Q} is such that $\limsup_{p \to \infty} C(p, \alpha) < \infty$, it is enough to assume $p < \infty$ first and conlude the case of $p$ infinite by letting $p \to \infty$. We will only prove \eqref{eqn:a_priori_poly} because the stability estimate \eqref{eqn:a_priori_poly_stab} comes from similar calculations.
		\medskip
		
		We start by assuming $p < \infty$. Applying~\eqref{eqn:coercivity} to the identity
		\begin{equation*}
			\frac{\d}{\d t}\bigg( \|\solp\|^p_{\spcp} \bigg)= p \int_{\R^d} \left\langle \partial_t \solp,\solp\right\rangle_{H^s_x} \|\solp\|_{H^s_x}^{p-2} \langle v \rangle^{p \alpha} \d v,
		\end{equation*}
		we obtain the differential inequality
		\begin{align*}
			\frac{\d}{\d t}\left( \frac{1}{p} \|\solp\|^p_{\spcp}\right) &+ \frac{\sigma_\BB}{\eps^2} \|\solp\|^p_{\spcpp} \\ \leq &\, \frac{1}{\eps} \int_{\R^d} \left\langle Q(h, h)+ Q(h, g),\solp \right\rangle_{H^s_x} \|\solp\|^{p - 2}_{H^s_x} \langle v \rangle^{p \alpha} \d v\\
			\leq &  \, \frac{C}{\eps} \bigg[ \|h\|_{\spcp} \| h\|_\spcpp + \|h\|_{\spcpp} \| g\|_\spcp + \|h\|_{\spcp} \| g\|_\spcpp \bigg] \|\solp\|^{p-1}_{\spcpp},
		\end{align*}
		where we used Lemma \ref{lem:estimate_Q}, from which the constant $C=C(p, \alpha)$ is from and satisfies $\limsup_{p \to \infty} C(p, \alpha) < \infty$. Multiplying both sides by $\exp(t p \sigma/\eps^2)$ (recall $\sigma \in (0, \sigma_\BB)$) and using the continuous inclusion $\Spcg{\beta} \subset \spcp$, we get
		\begin{align*}
			\frac{\d}{\d t}\left(\frac{1}{p} e^{\sigma p t /\eps^2} \|\solp\|_\spcp^p\right) +& \frac{\sigma_\BB - \sigma}{\eps^2} e^{\sigma p t / \eps^2} \|\solp\|_\spcpp^p \\
			& \leq \frac{C}{\eps} \left( \| h\|_\spcp  + \| g\|_\Spcg{\beta} \right) \left(e^{\sigma t /\eps^2} \|h\|_{\spcpp}\right) \left(e^{\sigma t /\eps^2} \|\solp\|_{\spcpp}\right)^{p-1}.
		\end{align*}
		Integrating and using Hölder's inequality with the exponents $\left(\infty, p, \frac{p}{p-1}\right)$, one gets
		\begin{align*}
			\tspcpN{\solp}^p \leq  \frac{C}{\eps} \bigg(\tspcpN{h} + \|g\|_{L^\infty_t \Spcg{\beta}}&\bigg)
			\times \sup_{0 \leq t \leq T}  \left\{\left(\int_0^t \left(e^{\sigma t'/\eps^2} \|h(t')\|_{\spcpp}\right)^p \d t'\right)^{1/p} \right. \\
			& \times \left. \left(\int_{0}^{t} \left(e^{\sigma t' /\eps^2} \|\solp(t')\|_{\spcp}\right)^p \d t'\right)^{\frac{p-1}{p}} \right\},
		\end{align*}
		from which we deduce
		\begin{align*}
			\tspcpN{\solp}^p &\leq \frac{C}{\eps} \bigg(\tspcpN{h} + \|g\|_{L^\infty_t \Spcg{\beta}}\bigg) \left(\eps^2 \tspcpN{h}^p \right)^{1/p} \left( \eps^2 \tspcpN{\solp}^p \right)^{\frac{p-1}{p}} \\
			& \leq C \eps \bigg(\tspcpN{h} + \|g\|_{L^\infty_t \Spcg{\beta}}\bigg) \tspcpN{h} \tspcpN{\solp}^{p-1}.
		\end{align*}
		This concludes the proof.
	\end{proof}
	
	\section{Bounds and asymptotics of the semigroup}
	\resetcounter
	
	\label{scn:asymptotics_semigroup}
	
	In this section, we establish integrability properties and asymptotics on the semigroup~$U^\eps$, necessary to prove the vanishing of the coupling term $\eps^{-2} U^\eps * \AA h^\eps$ in~$L^1_t \Spcg{\beta} + L^\infty_t \Spcg{\beta}$ (Lemmas \ref{lem:est_convo_1}-\ref{lem:est_convo_infty}) and the vanishing of the source term $\SS^\eps[h^\eps]$ from~\eqref{eqn:coupled_diff} (Lemma \ref{lem:cv_source_NSFI}). This will require to generalize some estimates proved in \cite{IGT20, BU91} in $\Spcg{\beta}$ to the larger space $\Spcp{p}{\alpha}$. To do so, we draw inspiration from the factorization techniques used in \cite[Section 5]{BU91}, using this time the splitting of the linearized operator~$\LL$ recalled in Section \ref{scn:splitting_op}. For the sake of completeness, we prove some results already present in \cite{BU91, IGT20}.
	
	\subsection{The eigenprojectors and partial semigroups on the Gaussian space}
	Let us first recall the spectral study originally led in \cite{EP75} for hard cutoff potentials in the space~$L^2\left(M^{-1/2}\right)$. We also mention the founding paper \cite{Nico} which has initiated the study of the spectrum of the linearized Boltzmann operator, and such works as \cite{G21, YY16, R13} in which these results were partially generalized.

	\medskip
	
	For a family of operator $(T(\xi) )_{\xi \in \R^d}$ acting on the $v$-variable, we define the Fourier multipliers~$T(D)$ acting on functions $u=u(x, v)$ by
	$$\widehat{T(D) u}(\xi, v) := \bigg(T(\xi) \widehat{u}(\xi , \cdot) \bigg)(v),$$
	where $u \mapsto \widehat{u}$ represents the Fourier transform with respect to the variable $x \in \Omega$. Note that such operators commute with $D = -i \frac{\d}{\d x}$ and any Fourier multiplier $f(D)$.
	\label{scn:spectral_splitting}

	\medskip
	\noindent
	\textbf{Spectral decomposition and expansions in Fourier space.} According to \cite[Theorems 1 and 2]{G21}, there exist $\err > 0$ (which can be assumed as small as necessary) and $\ah=\ah(\err) \in (0, \nu_0)$, where $\nu_0$ was presented in Section \ref{scn:boltzmann}, a family of projectors $\left( \PP^{(\ell)}_{\flat, j}(\xi) \right)_{\xi \in \R^d}$ uniformly bounded in~$\BBB\left( L^2_v\left(M^{-1/2}\right) \right)$, complex numbers $(\lambda_j(\xi))_{\xi \in \R^d}$, and closed operators~$\left( \LL_\sharp(\xi) \right)_{\xi \in \R^d}$ in~$\CCC\left( L^2_v\left( M^{-1/2} \right) \right)$, with $j=-1, \dots, 2$ and $\ell =0, 1, 2$, such that the following spectral decomposition holds in Fourier space for any $\xi \in \R^d$:
	\begin{gather}
		\label{eqn:eigen_relation}
		\LL + i v \cdot \xi = \sum_{j = -1}^2 \lambda_j(\xi) \PP_{\flat, j}(\xi) + \LL_\sharp(\xi),\\
		\label{eqn:ortho_proj}
		\PP_{\flat, j} \LL_\sharp = \LL_\sharp \PP_{\flat, j} =  \PP_{\flat, j} \PP_{\flat, k} = 0, ~ j \neq k,
	\end{gather}
	and the operator $\LL_\sharp(\xi) + \ah$ generates a bounded $\CC^0$-semigroup, uniformly in $\xi$. These eigenprojectors and eigenvalues expand around $\xi = 0$:
	\begin{gather}
		\label{eqn:eigen_expansion}
		\lambda_j(\xi) = \lambda_j^{(1)} |\xi|+ \lambda_j^{(2)} |\xi|^2 + O(|\xi|^3), ~
		\lambda_j^{(1)} \in i \R, ~ \lambda_j^{(2)} < 0,\\
		\label{eqn:proj_expansion}
		\PP_{\flat, j}(\xi) = \PP_{\flat, j}^{(0)}\left(\widetilde{\xi}\right) + |\xi| \PP_{\flat, j}^{(1)}\left(\xi\right), ~ \widetilde{\xi} := \frac{\xi}{|\xi|},
	\end{gather}
	and the zeroth order eigenprojectors sum to the $L^2_v(M^{-1/2})$ orthogonal projection on the null space of $\LL$:
	\begin{equation}
		\label{eqn:sum_proj}
		\Pi = \sum_{j=-1}^2 \PP_{\flat, j}^{(0)}\left( \widetilde{\xi} \right).
	\end{equation}
	In \cite{G21}, the projectors $\PP_{\flat, j}(\xi)$ and eigenvalues $\lambda_j(\xi)$ are actually defined for small frequencies $|\xi| \leq \err$, we assume here for simplicity that they are defined for any $\xi \in \R^d$ but vanish for $|\xi| > \err$, and we denote $\chi$ the characteristic function of $\{|\xi| \leq \err\}$. 
	
	\bigskip
	\noindent
	\textbf{Transposition to the physical space.} Let us denote by $\spcgg$ the following weighted Sobolev space:
	\begin{equation}
		\label{eqn:def_spcgg}
		\spcgg := L^2_v H^s_x \left(M^{-1/2}\right) = H^s_x L^2_v \left(M^{-1/2}\right),
	\end{equation}
	endowed with the norm
	\begin{align*}
		\left\| f \right\|_{\spcgg}^2 = \int_{\R^d \times \R^d} \left|\widehat{f}(\xi, v)\right|^2 \langle \xi \rangle^{2 s} M^{-1}(v) \d v \d \xi, \text{ if } \Omega = \R^d,\\
		\left\| f \right\|_{\spcgg}^2 = \sum_{\xi \in \Z ^d}  \langle \xi \rangle^{2 s} \int_{\R^d} \left|\widehat{f}(\xi, v)\right|^2 M^{-1}(v) \d v, \text{ if } \Omega = \T^d.
	\end{align*}
	Following the previous spectral decomposition in Fourier space, we define for the scaled linearized operator $\eps^{-2}\left(\LL + i \eps v \cdot \xi\right)$ the \textit{approximate hydrodynamic projector} and associated partial semigroup:
	\begin{gather}
		\label{eqn:sum_hyproj}
		\hyproj := \sum_{j=-1}^{2} \PP_{\flat, j}(\eps D),\\
		\label{eqn:semigroup_branch_def}
		U^\eps_\flat:= U^\eps \hyproj = \hyproj U^\eps= \sum_{j = -1}^2  \exp \left(t \lambda_j(\eps D)/\eps^2\right) \PP_{\flat, j}(\eps D).
	\end{gather}
	Since the operator $\LL$ is non-positive on $L^2_v\left(M^{-1/2}\right)$ (see \cite[Proposition 2.11]{UY06}) and $v \cdot \nabla_x$ is skew-adjoint in $H^s_x$ for any $v \in \R^d$, the scaled linearized operator~$\eps^{-2}\left(\LL + \eps v \cdot \nabla_x\right)$ is non-positive and thus $U^\eps$ is a contraction semigroup on $\spcgg$. We deduce by the uniform bounds on $\PP_{\flat, j}^{(\ell)}(\xi)$ that for some $C > 0$
	\begin{gather}
		\notag
		\left\| \hyproj \right\|_{ \BBB(\spcgg) } \leq C,\\
		\label{eqn:L2_SG_hydro}
		\left\| U_\flat^\eps(t) \right\|_{ \BBB(\spcgg) } \leq C.
	\end{gather}
	We also define the complementary projector and partial semigroup:
	\begin{gather*}
		\hyprojo := \textnormal{Id} - \hyproj,\\
		U^\eps_\sharp := U^\eps \hyprojo = \hyprojo U^\eps  = \exp\left(t \LL_\sharp(\eps D) / \eps^2\right).
	\end{gather*}
	By the boundedness of $\hyproj$ and the fact that $\LL_\sharp(D) + \ah$ generates a bounded $\CC^0$ semigroup on $\spcgg$, we may also assume that the constant $C$ is such that
	\begin{gather}
		\notag
		\left\| \hyprojo \right\|_{ \BBB(\spcgg) } \leq C,\\
		\label{eqn:L2_dec_SG_hydroo}
		\left\| U^\eps_\sharp(t) \right\|_{\BBB(\spcgg)} \leq C e^{-\ah t/\eps^2}.
	\end{gather}
	
	\medskip
	\noindent
	\textbf{Asymptotic behavior of $\boldsymbol{U^\eps_\flat}$.}
	The coefficients in the expansion of the scaled eigenvalues $\lambda_j(\xi)$ and the corresponding eigenmodes indicate (see \cite{EP75, G21}, \cite[Proposition A.3]{IGT20} and \cite[Remark 2.2.12]{UY06}) that, for $j=0, 2$, the eigenvalues correspond to the diffusion terms in \eqref{eqn:INSF}, and the projectors to the corresponding subspaces of macroscopic distributions $\bigg(\rho(x) + u(x) \cdot v + \frac{1}{d}\left(|v|^2 - d\right)\bigg)M(v)$:
	\begin{gather*}
		\eps^{-2} \lambda_0(\eps D) \approx - \lambda_0^{(2)} \Delta_x = \kappa \Delta_x,\\
		\PP_{\flat, 0}\left( \eps D \right) \approx \PP_{\flat, 0}^{(0)}\left( \widetilde{D} \right) = \text{$\spcgg$-orthogonal projection on } \rho + \theta = 0,\\
		\eps^{-2} \lambda_2(\eps D) \approx - \lambda_2^{(2)} \Delta_x =\mu \Delta_x,\\
		\PP_{\flat, 2}\left( \eps D \right) \approx \PP_{\flat, 2}^{(0)}\left( \widetilde{D} \right) = \text{$\spcgg$-orthogonal projection on } \nabla_x \cdot u = 0,
	\end{gather*}
	where we recall that $\mu$ and $\kappa$ are respectively the kinematic viscosity and thermal conductivity of the fluid represented by $M$, and $\widetilde{D}$ is the pseudodiferential operator associated witjh the symbol $\widetilde{\xi} = \xi / |\xi|$. This means that $\PP_{\flat, 0}^{(0)}\left( \widetilde{D} \right) + \PP_{\flat, 2}^{(0)}\left( \widetilde{D} \right)$ is the  projector on well-prepared distributions defined in \eqref{eqn:well_prepared_initial data}, and also that we have the following expression for the semigroup of the incompressible Navier-Stokes-Fourier system in its kinetic formulation~\eqref{eqn:KINSF}, presented in \cite{BU91}:
	\begin{gather}
		\label{eqn:def_U_0}
		U^0(t) := e^{\kappa t \Delta_x} \PP_{\flat, 0}^{(0)}\left( \widetilde{D} \right) + e^{\mu t \Delta_x} \PP_{\flat, 2}^{(0)}\left( \widetilde{D} \right),\\
		\label{eqn:def_Psi_0}
		\Psi^0(t)(f, g) := \int_{0}^{t} U^0\left(t - t'\right) \left(\PP_{\flat, 0}^{(1)} \left(\widetilde{D}\right) + \PP_{\flat, 2}^{(1)} \left(\widetilde{D}\right) \right) Q\left( f(t'), g(t') \right) \d t'.
	\end{gather}
	Similarly, for $j =\pm 1$, denoting $c$ the speed of sound in the gas represented by $M$, we have the following asymptotics:
	\begin{align*}
		\eps^{-2} \lambda_{\pm 1}(\eps D) \approx& \pm i \frac{c}{\eps} |D|  + \kappa \Delta_x,\\
		\PP_{\flat, 1}^{(0)}\big( \eps D \big) + \PP_{\flat, -1}^{(0)}\left( \eps D \right) \approx&
		\PP_{\flat, 1}^{(0)}\big( \widetilde{D} \big) + \PP_{\flat, -1}^{(0)}\left( \widetilde{D} \right)\\
		=& \text{$\spcgg$-orthogonal projection on}\\
		&\text{ill-prepared macro. distributions.}
	\end{align*}
	We define the semigroup $U^\eps_\disp$ corresponding to these acoustic waves through
	\begin{equation}
		\label{eqn:def_U_eps_disp}
		U^\eps_\disp(t)  e^{-t \kappa \Delta_x} := e^{i t \frac{c}{\eps} |D|} \PP_{\flat, 1}^{(0)}\left( \widetilde{D} \right) + e^{-i t \frac{c}{\eps} |D|} \PP_{\flat, -1}^{(0)}\left( \widetilde{D} \right),
	\end{equation}
	so as to highlight the presence of the wave operator. Using the notations of Section~\ref{scn:micro_macro_decomposition}, the orthogonality relation~\eqref{eqn:ortho_proj} implies that
	\begin{gather}
		\notag
		U^0 U^\eps_\disp = U^\eps_\disp U^0 = 0,\\
		\label{eqn:proj_hydro_semigroup}
		U^0 f = U^0 f_{\ini, \WP}, ~ U^\eps_\disp f = U^\eps_\disp f_{\ini, \IP}
	\end{gather}
	These hydrodynamic semigroups $U^0$ and $U^\eps_\disp$ were shown in \cite{BU91} to drive the dynamics of $U^\eps_\flat$, which is made explicit by this next lemma adapted from \cite[Lemma 3.5]{IGT20}
	\begin{lem}
		\label{lem:cv_U_eps_U_0}
		For any $\beta > d /2$ and $f \in \Spcg{\beta}$, there holds
		\begin{equation*}
			\lim\limits_{\eps \to 0} \left(\sup_{t \geq 0} \, \langle t \rangle^{1/2} \bigg\| \bigg(U^\eps_\flat(t) - U^0(t) - U^\eps_\disp(t) \bigg)\Pi f \bigg\|_{\Spcg{\beta}} \right)= 0.
		\end{equation*}
	\end{lem}
	\begin{proof}
		Recalling that $U^\eps_\flat = U^\eps - U^\eps_\sharp$ (our notations are consistent with \cite{IGT20}), this lemma comes simply from a density argument applied to \cite[Lemma 3.5, (3.12)-(3.13)]{IGT20} which we explain for clarity. This is possible because of the finite-dimensional (in $v$) nature of $\Pi f$:
		\begin{equation}
			\Pi f(x, v) = \bigg( \rho_f(x) + u_f(x) + \frac{1}{2}\left( |v|^2 - d \right) \theta_f(x)\bigg)M(v),
		\end{equation}
		where $\rho_f, u_f, \theta_f$ are $H^s_x$ functions and thus can be approximated by smooth functions.
	\end{proof}
	Let us also recall that the authors of \cite{BU91} proved that $\Psi^\eps(f, g) \rightarrow \Psi^0(f, g)$ in a weaker topology than that of $L^\infty_t \Spcg{\beta}$, and those of \cite{IGT20} proved it when $f=g=f^0$ by using the properties of the limit system \eqref{eqn:INSF}, we cite \cite[Lemma 4.1]{IGT20} here
	\begin{lem}
		For any $\beta > d/2 + 1$, there holds
		\begin{equation*}
			\lim\limits_{\eps \to 0} \bigg(\sup_{0 \leq t < T} \chi_\Omega(t) \left\| \Psi^\eps(f^0, f^0) - \Psi^0(f^0, f^0)\right\|_{ \Spcg{\beta} }\bigg) = 0,
		\end{equation*}
		where $\chi_\Omega(t) = \langle t \rangle^{1/4}$ if $\Omega = \R^2$, and $\chi_\Omega(t) = 1$ otherwise.
	\end{lem}
	Finally, the acoustic part was shown in \cite{IGT20} to be bounded, decay for large times, and vanish in a weaker topology when the domain $\Omega$ is the whole space, as explained by this lemma.
	\begin{lem}
		\label{lem:disp}
		For any $\beta > \frac{d}{2}$, there exists $C > 0$ such that
		\begin{gather}
			\label{eqn:Linf_bound_disp}
			\|U^\eps_\disp(t) f\|_{\Spcg{\beta}} \leq C \|f\|_{\Spcg{\beta}},\\
			\label{eqn:Linf_decay_disp}
			\langle t \rangle^{d/4} \|U^\eps_\disp(t) f\|_{{\Spcg{\beta}}} \leq C \left( \|f\|_{L^2_v L^1_x \left( M^{-1/2} \right) } + \|f\|_{{\Spcg{\beta}}} \right),
		\end{gather}
		furthermore, if $\Omega = \R^d$ and $f$ is a Schwartz function, there holds for some $C_f > 0$
		\begin{equation}
			\label{eqn:Linf_statphase_disp_bis}
			\|U^\eps_\disp(t) f \|_{L^\infty_v W^{s, \infty}_x \left( \langle v \rangle^{\beta} M^{-1/2} \right)}
			\leq C_f \left(\frac{\eps}{t}\right)^{\frac{d-1}{2}}.
		\end{equation}
	\end{lem}

	Note that the estimate \eqref{eqn:Linf_bound_disp} comes from \eqref{eqn:def_U_eps_disp} and the fact that $\PP_{\flat, \pm 1}^{(0)} \left( \widetilde{D} \right)f$ is a linear combination of functions of the form
	\begin{equation*}
		\bigg(\int_{\R^d } \langle f(v_*, \cdot), g \rangle_{H^s_x} P_1(v_*) M(v_*) \d v_*\bigg) P_2 M
	\end{equation*}
	for some polynomials $P_1, P_2$ and $g \in H^s_x$.
	The estimate \eqref{eqn:Linf_decay_disp} is \cite[(3.25)]{IGT20}, and~\eqref{eqn:Linf_statphase_disp_bis} is a direct consequence of the last estimate of \cite[p. 587]{IGT20}.

	\begin{figure}
		\centering
		
		\caption{Spectrum of $\frac{1}{\eps^2}\left( \LL + \eps v \cdot \nabla_x \right)$}
		
		
		\tikzset{
			pattern size/.store in=\mcSize, 
			pattern size = 5pt,
			pattern thickness/.store in=\mcThickness, 
			pattern thickness = 0.3pt,
			pattern radius/.store in=\mcRadius, 
			pattern radius = 1pt}
		\makeatletter
		\pgfutil@ifundefined{pgf@pattern@name@_4l82bz2cu lines}{
			\pgfdeclarepatternformonly[\mcThickness,\mcSize]{_4l82bz2cu}
			{\pgfqpoint{0pt}{0pt}}
			{\pgfpoint{\mcSize+\mcThickness}{\mcSize+\mcThickness}}
			{\pgfpoint{\mcSize}{\mcSize}}
			{\pgfsetcolor{\tikz@pattern@color}
				\pgfsetlinewidth{\mcThickness}
				\pgfpathmoveto{\pgfpointorigin}
				\pgfpathlineto{\pgfpoint{\mcSize}{0}}
				\pgfusepath{stroke}}}
		\makeatother
		\tikzset{every picture/.style={line width=0.75pt}} 
		
		\begin{tikzpicture}[x=0.75pt,y=0.75pt,yscale=-1,xscale=1]
			
			\fill[pattern=north west lines, pattern color=black] (330,380) rectangle (240,180);
			\draw  [draw opacity=0][fill={rgb, 255:red, 64; green, 64; blue, 64 }  ,fill opacity=1 ] (150,180) -- (240,180) -- (240,380) -- (150,380) -- cycle ;
			\draw [color={rgb, 255:red, 0; green, 0; blue, 255 }  ,draw opacity=1 ][line width=1.5]    (330,280) -- (440,280) ;
			\draw    (440,280) -- (480,280) ;
			\draw [color={rgb, 255:red, 0; green, 0; blue, 0 }  ,draw opacity=1 ][pattern=_4l82bz2cu,pattern size=6pt,pattern thickness=0.75pt,pattern radius=0pt, pattern color={rgb, 255:red, 0; green, 255; blue, 0}]   (330,170) -- (330,390) ;
			\draw  [draw opacity=0] (330.38,350) .. controls (390,350) and (440,320) .. (440,280) .. controls (440,240) and (390,210) .. (330,210) -- (330,280) -- cycle ; \draw  [color={rgb, 255:red, 255; green, 0; blue, 0 }  ,draw opacity=1 ] (330,350) .. controls (390,350) and (440,320) .. (440,280) .. controls (440,240) and (390,210) .. (330,210) ;
			\draw  [draw opacity=0][dash pattern={on 4.5pt off 4.5pt}] (380,380) .. controls (415,360) and (440,325) .. (440,280) .. controls (440,235) and (415,200) .. (380,180) -- (330,280) -- cycle ; \draw  [color={rgb, 255:red, 255; green, 0; blue, 0 }  ,draw opacity=1 ][dash pattern={on 4.5pt off 4.5pt}] (380,380) .. controls (415,360) and (440,320) .. (440,280) .. controls (440,240) and (415.63,199.7) .. (380.12,181.63) ;
			\draw  [draw opacity=0][dash pattern={on 0.84pt off 2.51pt}] (416.28,377) .. controls (431.07,350.85) and (440,315) .. (440,280) .. controls (440,240) and (430,210) .. (416.52,183.42) -- (340,280) -- cycle ; \draw  [color={rgb, 255:red, 255; green, 0; blue, 0 }  ,draw opacity=1 ][dash pattern={on 0.84pt off 2.51pt}] (416.28,377) .. controls (431.07,350.85) and (440,316.99) .. (440,280) .. controls (440,243.21) and (431.17,209.52) .. (416.52,183.42) ;
			\draw    (440,170) -- (440,390) ;
			\draw [color={rgb, 255:red, 0; green, 0; blue, 0 }  ,draw opacity=1 ][line width=0.75]    (240,170) -- (240,390) ;
			
			\draw (451,292.4) node [anchor=north west][inner sep=0.75pt]    {$0$};
			\draw (430,390) node [anchor=north west][inner sep=0.75pt]    {$i\mathbb{R}$};
			\draw (463.25,188.89) node [anchor=north west][inner sep=0.75pt]  [color={rgb, 255:red, 255; green, 0; blue, 0 }  ,opacity=1 ]  {$\lambda _{+1}( \varepsilon \xi ) /\varepsilon ^{2}$};
			\draw (463.25,343.46) node [anchor=north west][inner sep=0.75pt]  [color={rgb, 255:red, 255; green, 0; blue, 0 }  ,opacity=1 ]  {$\lambda _{-1}( \varepsilon \xi ) /\varepsilon ^{2}$};
			\draw (490,230) node [anchor=north west][inner sep=0.75pt]  [color={rgb, 255:red, 0; green, 0; blue, 255 }  ,opacity=1 ]  {$\lambda _{0}( \varepsilon \xi ) /\varepsilon ^{2}$};
			\draw (310,390) node [anchor=north west][inner sep=0.75pt]    {$-\frac{\mathbf{a}}{\varepsilon ^{2}} +i\mathbb{R}$};
			\draw (200,390) node [anchor=north west][inner sep=0.75pt]  [color={rgb, 255:red, 0; green, 0; blue, 0 }  ,opacity=1 ]  {$-\frac{\nu _{0}}{\varepsilon ^{2}} +i\mathbb{R}$};
			\draw (490,290) node [anchor=north west][inner sep=0.75pt]  [color={rgb, 255:red, 0; green, 0; blue, 255 }  ,opacity=1 ]  {$\lambda _{2}( \varepsilon \xi ) /\varepsilon ^{2}$};
			
		\end{tikzpicture}
		
		The blue part corresponds to \eqref{eqn:INSF}. The red part, which gets closer to the imaginary axis as $\eps \to 0$, corresponds to the acoustic waves. The hatched part may contain spectral values, and the solid part is included in the spectrum.
	\end{figure}
	
	\subsection{Spectral properties on the polynomial spaces}
	We now extend some of the properties of $U^\eps_\flat$ proved on $\BBB(\spcgg)$ in \cite{IGT20} to $\BBB\left( \Spcp{p}{\alpha} ; \Spcg{\beta} \right)$. This will be made possible using factorization techniques with the Gualdani-Mischler-Mouhot decomposition~$\LL = \AA + \BB$ and \textit{Grad's decomposition} (see \cite{G63, UY06}):
	\begin{equation}
		\label{eqn:grad_decomposition}
		\LL = -\nu + K
	\end{equation}
	where $\nu$ was defined by \eqref{eqn:def_coll_freq}, and $K$ satisfies the following regularization property:
	\begin{equation}
		\label{eqn:reg_K}
		\forall \beta \geq 0, ~ K \in \BBB\left( \spcgg ; \Spcg{0} \right) \cap \BBB\left( \Spcg{\beta} ; \Spcg{\beta + 1} \right).
	\end{equation}
	As a first step, we will study operators $\RR^\eps_j$ and $R^\eps_j$ (Lemmas \ref{lem:regularizer} and \ref{lem:regularizer_Linf}) which appear in factorization formulas (\eqref{eqn:fact_1} and \eqref{eqn:fact_2}) for $\PP_{\flat, j}(\eps D)$, this will allow to prove boundedness properties (Lemma \ref{lem:bound_hyproj_L2} and \ref{lem:bound_hyproj_Linf}) for the coefficients in the expansion (coming from \eqref{eqn:eigen_expansion})
	$$\PP_{\flat, j}(\eps D) = \PP_{\flat, j}^{(0)}\left(\widetilde{D}\right) + \eps |D| \PP^{(1)}_{\flat, j}(\eps D).$$
	These properties will be used to prove bounds and asymptotics on $U^\eps_\flat$ and $U^\eps_\sharp$ (Lemmas \ref{lem:est_U_opt}-\ref{lem:cross_partial_semigroup}).

	\medskip
	We recall that $\alpha_\BB, \AA, \BB, \BB^\eps$ were defined in Section \ref{scn:splitting_op}, and that
	\begin{equation*}
		\BB^\eps = \frac{1}{\eps^2} \left( \BB + \eps \cdot \nabla_x\right).
	\end{equation*}
	
	\begin{lem}
		\label{lem:regularizer}
		For any $p\in [1, \infty]$, $\alpha > \alpha_\BB$, $j=-1,\dots, 2$, the following operator is well-defined and bounded uniformly in $\eps \in (0, 1)$:
		\begin{equation*}
			\RR_j^\eps := \frac{1}{\eps^2} \AA \left(\eps^{-2} \lambda_j(\eps D) - \BB^\eps\right)^{-1} \chi(\eps D) \in \BBB\left( \Spcp{p}{\alpha} ; \spcgg \right),
		\end{equation*}
		where the inverse is to be understood on the range $\chi(\eps D) \Spcp{p}{\alpha}$. It expands as
		\begin{equation*}
			\RR_j^\eps = -\AA \BB^{-1} + \eps |D| \RR^{\eps, 1}_j,
		\end{equation*}
		with $\RR_j^{\eps, 1} \in \BBB(\Spcp{p}{\alpha} ; \spcgg)$ uniformly in $\eps \in (0, 1)$, and every term commutes with $D$.
	\end{lem}

	\newcommand{\Spcpeps}{\spcp^{p, \alpha}_{\err / \eps}}
	\begin{proof}
		In this proof, we drop the subscript $j$ and denote
		\begin{gather*}
			\Spcpeps := \chi(\eps D) \Spcp{p}{\alpha} = \bigg\{ f \in \Spcp{p}{\alpha} ~ : ~ \widehat{f}(\xi)=0 \text{ if } \eps |\xi| > \err\bigg\},\\
			T^\eps :=\eps v \cdot \nabla_x - \lambda(\eps D),
		\end{gather*}
		so that the regularizing operator rewrites
		\begin{equation}
			\label{eqn:cara_R_eps}
			\RR^\eps = - \AA \left( \BB + T^\eps \right)^{-1} \chi(\eps D).
		\end{equation}
		In the two first steps, we will prove the following properties:
		\begin{gather}
			\label{eqn:H1}
			\BB \nu^{-1} \text{ is invertible on $\Spcpeps$ and its inverse is bounded by $c > 0$},\\
			\label{eqn:H2}
			\left\| T^\eps \nu^{-1} \right\|_{\BBB\left( \Spcpeps \right)} \leq \frac{1}{2 c},\\
			\label{eqn:H3}
			T^\eps = \eps |D| T^{\eps, 1} \text{ and $\left\|T^{\eps, 1} \nu^{-1}\right\|_{\BBB\left( \Spcpeps \right)} \leq c$},
		\end{gather}
		and then use them in a third step to prove the lemma.
		
		\step{1}{Proof of \eqref{eqn:H1}}
		According to \cite[(4.33)]{GMM17}, $\BB - \nu$ is $\nu$-bounded with relative bound less that 1 (see Lemma \ref{lem:diss_B_eps}), thus, $$\BB \nu^{-1} = \Id + (\BB - \nu)\nu^{-1},$$
		which implies \eqref{eqn:H1} because $\left\| \left( \BB - \nu \right) \nu^{-1} \right\|_{\BBB\left( \Spcpeps \right)}  < 1$.
		
		\step{2}{Proof of \eqref{eqn:H2} and \eqref{eqn:H3}}
		By \eqref{eqn:eigen_expansion}, $T^\eps \nu^{-1}$ is a Fourier multiplier in $x$ whose symbol is of the form
		\begin{align*}
			-i \eps\xi \cdot \frac{v}{\nu} - \frac{\lambda^{(1)}}{\nu} &\eps |\xi| + O\left(\frac{\eps^2 |\xi|^2}{\nu}\right)\\
			&= \eps |\xi| \left( \widetilde{\xi} \cdot \frac{v}{\nu} - \frac{\lambda^{(1)} }{\nu} + O\left( \eps |\xi| \right) \right),
		\end{align*}
		so there holds indeed
		$$T^\eps = \eps |D| T^{\eps, 1}.$$
		Furthermore, on $\Spcpeps$, we have by definition $\eps |\xi| \leq \err$, it is then clear that the operators $T^\eps \nu^{-1}$ and $T^{\eps, 1} \nu^{-1}$ are bounded on $\Spcpeps$ uniformly in $\eps \in (0, 1)$. Finally, up to  a reduction of $\err$, there also holds $\left\|T^\eps \nu^{-1}\right\|_{\BBB\left( \Spcpeps \right)} \leq 1/2c$.
		
		\step{3}{Proof of the lemma}
		Because of assumptions \eqref{eqn:H1} and \eqref{eqn:H2}, the following operator is well-defined on $\Spcpeps$ and bounded by $1/2 \nu_0$:
		\begin{equation*}
			\left( \BB + T^\eps \right)^{-1} = \nu^{-1} \left( \BB \nu^{-1} + T^\eps \nu^{-1}\right)^{-1}
		\end{equation*}
		Combined with \eqref{eqn:cara_R_eps} and the regularization property \eqref{eqn:reg_AA}, this yields the boundedness property $\RR^\eps \in \BBB\left( \Spcp{p}{\alpha} ; \spcgg \right)$ uniformly in $\eps \in (0, 1)$. Let us now tun to its expansion; there holds on $\Spcpeps$:
		\begin{align*}
			\left(\BB + T^\eps\right)^{-1} &= \BB^{-1} + \BB^{-1} \left(T^{\eps} \nu^{-1}\right) \left( \BB \nu^{-1} + T^\eps \nu^{-1} \right)^{-1}\\
			&= \BB^{-1} + \eps |D| \BB^{-1} \left(T^{\eps, 1} \nu^{-1}\right) \left( \BB \nu^{-1} + T^\eps \nu^{-1} \right)^{-1}\\
			&=: \BB^{-1} + \eps |D| S^\eps,
		\end{align*}
		where $S^\eps \in \BBB( \Spcpeps )$ uniformly in $\eps \in (0, 1)$. We deduce that
		\begin{align*}
			\RR^\eps &= -\AA \left( \BB + T^\eps \right)^{-1} \chi(\eps D)\\
			&= -\AA \BB \chi(\eps D) + \eps |D| S^\eps \chi(\eps D)\\
			&= -\AA \BB + \eps |D| \RR^{\eps, 1},
		\end{align*}
		where we let
		\begin{equation*}
			\RR^{\eps, 1} := S^\eps \chi(\eps D) + \left(\eps |D|\right)^{-1} \left( \Id - \chi(\eps D) \right) \AA \BB^{-1},
		\end{equation*}
		which is a bounded operator in virtue of the inequality $1 - \chi(\eps \xi) \leq \eps |\xi| / \err$.

	\end{proof}

	\begin{lem}
		\label{lem:bound_hyproj_L2}
		For any $p\in [1, \infty]$ and $\alpha > \alpha_\BB$, there exists $C > 0$ such that
		\begin{gather}
			\label{eqn:bound_hyproj_branch}
			\left\| \PP_{\flat, j} (\eps D) f\right\|_{\spcgg}  \leq C \|f\|_{\Spcp{p}{\alpha}},\\
			\label{eqn:bound_hyproj_derivative}
			 \left\| \PP_{\flat, j}^{(1)} (\eps D) f\right\|_{\spcgg} \leq C \|f\|_{\Spcp{p}{\alpha}}.
		\end{gather}
		where we recall that $\spcgg = L^2_v H^s_x \left(M^{-1/2}\right)$.
	\end{lem}
	
	\begin{proof}
		From the relations \eqref{eqn:eigen_relation}-\eqref{eqn:ortho_proj} and the splitting $\LL^\eps = \eps^{-2} \AA + \BB^\eps$, there holds
		\begin{align*}
			\PP_{\flat, j}(\eps D)  \LL^\eps = \PP_{\flat, j}(\eps D) \eps^{-2} \lambda_j(\eps D)= \PP_{\flat, j}(\eps D) \left(\eps^{-2} \AA + \BB^\eps \right),
		\end{align*}
		which can be rearranged as
		\begin{equation*}
			\PP_{\flat, j}(\eps D) \left( \eps^{-2} \lambda_j(\eps D) - \BB^\eps \right) =  \eps^{-2}  \PP_{\flat, j}(\eps D) \AA,
		\end{equation*}
		and thus the following factorization formula holds:
		\begin{equation}
				\label{eqn:fact_1}
				\PP_{\flat, j}(\eps D) = \PP_{\flat, j}(\eps D) \RR_j^\eps,
		\end{equation}
		where $\RR_j^\eps$ is that of Lemma \ref{lem:regularizer}.
		Thanks to the boundeness of $\RR_j^\eps$ proved in Lemma~\ref{lem:regularizer}, and the boundedness~$\PP_{\flat, j}(\eps D) \in \BBB(\spcgg)$ recalled in the beginning of the section, we deduce that, uniformly in $\eps \in (0, 1)$, there holds $\PP_{\flat, j} (\eps D) \in \BBB\left( \Spcp{p}{\alpha} ; \spcgg \right)$.
		
		Furthermore, injecting the expansion from Lemma \ref{lem:regularizer} and the one of $\PP_{\flat, j}(\eps D)$ in~\eqref{eqn:fact_1}, we get
		\begin{align*}
			\PP_{\flat, j}(\eps D) &= - \PP_{\flat, j}^{(0)} \left(\widetilde{D}\right) \AA \BB^{-1} \\
			&+ \eps |D|\bigg( \PP_{\flat, j}^{(1)}\left(\eps D\right) \RR_j +  \PP_{\flat, j}^{(0)}\left(\widetilde{D}\right) \RR_j^{\eps, 1} \bigg),
		\end{align*}
		and thus, identifying the first order coefficient, we have
		\begin{align*}
			\PP_{\flat, j}^{(1)}(\eps D) = \PP_{\flat, j}^{(1)}\left(\eps D\right) \RR_j^\eps +  \PP_{\flat, j}^{(0)}\left(\widetilde{D}\right) \RR_j^{\eps, 1}.
		\end{align*}
		Since $\PP_{\flat, j}^{(1)}(\eps D) \in \BBB(\spcgg)$ uniformly in $\eps$, we conclude that $\PP_{\flat, j}^{(1)}(\eps D) \in \BBB(\Spcp{p}{\alpha} ; \spcgg)$ uniformly in $\eps \in (0, 1)$ thanks to Lemma \ref{lem:regularizer}.
		
	\end{proof}
	
	\begin{lem}
		\label{lem:regularizer_Linf}
		For any $j=-1,\dots, 2$, the following operator is well-defined:
		\begin{equation*}
			R_j^\eps := K \left( \lambda_j(\eps D) - \nu + \eps v \cdot \nabla_x\right)^{-1} \chi(\eps D),
		\end{equation*}
		furthermore, it expands as
		\begin{equation*}
			R_j^\eps = -K \nu^{-1} + \eps |D| R^{\eps, 1}_j,
		\end{equation*}
		where each term commutes with $D$ and for any $\beta \geq 0$, we have,
		\begin{equation*}
			R_j^\eps, R_j^{\eps, 1} \in \BBB\left(\spcgg ; \spcg^0\right) \cap \BBB\left( \Spcg{\beta} ; \Spcg{\beta+1} \right),
		\end{equation*}
		uniformly in~$\eps \in (0, 1)$.
	\end{lem}

	\begin{proof}
		This result can be proved in the same way as was Lemma \ref{lem:regularizer}, but using Grad's decomposition \eqref{eqn:grad_decomposition} and the role of $\AA$ (resp. $\BB$) is replaced by $K$ (resp. $\nu$).
		
	\end{proof}

	\begin{lem}
		\label{lem:bound_hyproj_Linf}
		For any $p\in [1, \infty]$, $\beta \geq 0$, $\alpha > \alpha_\BB$, there exists $C > 0$ such that
		\begin{gather}
			\label{eqn:bound_hyproj_branch_Linf}
			\left\| \PP_{\flat, j} (\eps D) f\right\|_{\Spcg{\beta}}  \leq C \|f\|_{\Spcp{p}{\alpha}},\\
			\label{eqn:bound_hyproj_derivative_Linf}
			\left\| \PP_{\flat, j}^{(1)} (\eps D) f\right\|_{\Spcg{\beta}} \leq C \|f\|_{\Spcp{p}{\alpha}},
		\end{gather}
		in particular, we have thanks to \eqref{eqn:sum_hyproj} that
		\begin{equation}
			\label{eqn:bound_hyproj}
			\| \hyproj f \|_{\Spcg{\beta}} \leq C \|f\|_\Spcp{p}{\alpha}.
		\end{equation}
	\end{lem}
		
%
		\begin{proof}
			First, note that by Lemma \ref{lem:regularizer_Linf}, there holds
			\begin{gather*}
				(R^\eps_j)^{1+\beta} = (-K \nu^{-1})^{1+\beta} + \eps |D| \widetilde{R}^\eps_j \in \BBB\left( \spcgg ; \Spcg{\beta} \right),\\
				\widetilde{R}^\eps_j \in \BBB\left( \spcgg ; \Spcg{\beta} \right),
			\end{gather*}
			uniformly in $\eps \in (0, 1)$. One shows a factorization formula in the spirit of \eqref{eqn:fact_1}, which iterated gives
			\begin{equation}
				\label{eqn:fact_2}
				\PP_{\flat, j}(\eps D) = \left(R_j^\eps\right)^{1+\beta} \PP_{\flat, j}(\eps D).
			\end{equation}
			We get \eqref{eqn:bound_hyproj_branch_Linf} by combining the bound of $(R^\eps_j)^{1+\beta}$ and \eqref{eqn:bound_hyproj_branch}. Let us now inject the expansion of $(R^\eps_j)^{1+\beta}$ and $\PP_{\flat, j}(\eps D)$ in the previous relation:
			\begin{align*}
				\PP_{\flat, j}(\eps D) =& \left( (-K \nu^{-1})^{1+\beta} + \eps |D| \widetilde{R}^\eps_j \right) \left( \PP_{\flat, j}^{(0)} \left( \widetilde{D} \right)  + \eps |D| \PP_{\flat, j}^{(1)} \left( \eps D \right)\right) \\
				=& (-K \nu^{-1})^{1+\beta}  \PP_{\flat, j}^{(0)} \left( \widetilde{D} \right) \\
				&+ \eps |D| \left( (-K \nu^{-1})^{1+\beta} \PP_{\flat, j}^{(1)} \left( \eps D \right) + \widetilde{R}^\eps_j \PP_{\flat, j} \left( \eps D \right)
				\right),
			\end{align*}
			which leads to the identification
			\begin{equation*}
				\PP_{\flat, j}^{(1)}\left( \eps D \right) = (-K \nu^{-1})^{1+\beta} \PP_{\flat, j}^{(1)} \left( \eps D \right) + \widetilde{R}^\eps_j \PP_{\flat, j} \left( \eps D \right),
			\end{equation*}
			where we have
			\begin{equation}
				\label{eqn:factorization_Linf}
				\left(R_j^\eps\right)^{1+\beta} \PP_{\flat, j}(\eps D) =  \PP_{\flat, j}(\eps D) \in \BBB\left(\spcgg ; \Spcg{\beta}\right).
			\end{equation}
			uniformly in $\eps \in (0, 1)$. Once again, we get \eqref{eqn:bound_hyproj_derivative_Linf} by combining the bound of $\widetilde{R}^\eps_j$ and~$(-K \nu^{-1})^{1+\beta}=\left(R_j^\eps\right)_{| \eps = 0}$ with \eqref{eqn:bound_hyproj_derivative}.
			
		\end{proof}

	\begin{cor}
		\label{lem:est_U_opt}
		For any $p\in [1, \infty]$, $\beta > d/2$ and $\alpha > \alpha_\BB$, there exists $C=C(p, \alpha, \beta)$ such that
		\begin{gather}
			\label{eqn:Linf_SG_hydro}
			\|U^\eps_\flat(t) f\|_{\Spcg{\beta}} \leq C \|f\|_{\Spcp{p}{\alpha}},\\
			\label{eqn:Linf_dec_SG_hydroo}
			\|U^\eps_\sharp(t) f\|_{\Spcg{\beta}} \leq C e^{-\ah t / \eps^2} \|f\|_{{\Spcg{\beta}}}.
		\end{gather}
	\end{cor}

	\begin{proof}
		Estimate \eqref{eqn:Linf_SG_hydro} comes from the combination of \eqref{eqn:L2_SG_hydro} with \eqref{eqn:bound_hyproj}.
		To get estimate \eqref{eqn:Linf_dec_SG_hydroo}, we apply the Duhamel formula \eqref{eqn:semigroup_fact} to Grad's decomposition \eqref{eqn:grad_decomposition}:
		\begin{equation*}
			U^\eps = S^\eps + \frac{1}{\eps^2} S^\eps * K U^\eps,
		\end{equation*}
		where we denoted
		\begin{align}
			\notag
			S^\eps(t) h (x, v) &:= \exp\left(\eps^{-2} t ( - \nu + \eps v \cdot \nabla_x)\right) h (x, v)\\
			\label{eqn:semigroup_transport}
			&= e^{- \frac{\nu(v) t}{\eps^2}} h(x - vt, v).
		\end{align}
		Compose from the right with $\hyprojo$ to get by the definition of $U^\eps_\sharp$
		\begin{equation*}
			U^\eps_\sharp = S^\eps \hyprojo + \frac{1}{\eps^2} S^\eps * K U^\eps_\sharp.
		\end{equation*}
		Note that in both cases $X= \spcgg$ or $E^\gamma$ (where $\gamma \geq 0$), we have thanks to \eqref{eqn:bound_col_freq}
		\begin{equation*}
			\left\| S^\eps(t) \right\|_{\BBB(X )} \lesssim \exp\left( -  \frac{\nu_0 t}{\eps^2}\right),
		\end{equation*}
		which implies in particular
		\begin{gather*}
			\frac{1}{\eps^2} \int_0^\infty \left\| S^\eps(t) \right\|_{\BBB(X )} \d t \lesssim \frac{1}{\eps^2} \int_0^\infty e^{- \nu_0 t / \eps^2} \d t \lesssim 1,\\
			\left\| S^\eps(t) \right\|_{\BBB(X )} \lesssim 1.
		\end{gather*}		
		Combined with the regularization property \eqref{eqn:reg_K}, we get	
		\begin{align*}
			\left\|U^\eps f \right\|_{L^\infty_t \spcg^{0}} &\lesssim \left\|S^\eps f \right\|_{L^\infty_t \spcg^{0}}  + \frac{1}{\eps^2} \int_0^t \left\| S^\eps(t-t') \right\|_{\BBB(\spcg^{0} )} \left\|U^\eps_\sharp(t') f \right\|_{L^\infty_t \spcgg} \d t' \\
			&\lesssim e^{-\nu_0 t / \eps^2}\left\|f \right\|_{\spcg^{0}} + \frac{1}{\eps^2} \int_{0}^{t} e^{-\nu_0 (t-t') / \eps^2} \left\|U^\eps_\sharp(t') f \right\|_{L^\infty_t \spcgg} \d t',
		\end{align*}
		and also for any $\gamma \geq 0$
		\begin{align*}
			\left\|U^\eps_\sharp(t) f \right\|_{L^\infty_t \spcg^{\gamma+1}} &\lesssim \left\|S^\eps(t) f \right\|_{\spcg^{\gamma+1}} + \frac{1}{\eps^2} \int_0^\infty \left\| S^\eps(t-t') \right\|_{\BBB(\spcg^{\gamma+1} )} \left\|U^\eps_\sharp(t') f \right\|_{\spcg^{\gamma}} \d t'\\
			& \lesssim e^{-\nu_0 t / \eps^2}\left\|f \right\|_{\spcg^{\gamma+1}} + \frac{1}{\eps^2} \int_{0}^{t} e^{-\nu_0 (t-t') / \eps^2} \left\|U^\eps_\sharp(t') f \right\|_{L^\infty_t \spcg^{\gamma}} \d t'
		\end{align*}
		Using the following relation valid for any $a, b > 0$,
		\begin{equation*}
			\int_{0}^{t} e^{-a (t-t')} e^{-b t'} \d t' \leq \frac{e^{- \min\{a, b\} t}}{|b-a|},
		\end{equation*}
		one gets by induction on $\gamma$, and using \eqref{eqn:L2_dec_SG_hydroo} as an initilization, that
		\begin{equation*}
			\left\| f \right\|_{\Spcg{\beta}} \lesssim e^{-\nu_0 t / \eps^2} \|f\|_{\Spcg{\beta}} + e^{-\min\{\ah, \nu_0\} t / \eps^2} \|f\|_{\spcgg}.
		\end{equation*}
		We conclude to \eqref{eqn:Linf_dec_SG_hydroo} using the continuous inclusion $\Spcg{\beta} \subset \spcgg$ as $\beta > d/2$.
	\end{proof}

	\begin{lem}
		\label{lem:est_u_pi_perp}
		For any $p\in [1, \infty]$, $\beta \geq 0$, $\alpha > \alpha_\BB$, there exists $C=C(p, \alpha, \beta)$ such that for any $f \in \Spcp{p}{\alpha}$
		\begin{equation*}
			\|U^\eps_\flat(t) \left(\textnormal{Id} - \Pi\right) f\|_{\Spcg{\beta}}  \leq C \min\left\{\|f\|_\Spcp{p}{\alpha},  \frac{\eps}{t^{1/2}} \|f\|_\Spcp{p}{\alpha} , \eps \|\nabla_x f\|_\Spcp{p}{\alpha} \right\}.
		\end{equation*}
	\end{lem}

	\begin{proof}
		First, recall that we have from \eqref{eqn:eigen_expansion}, and from \eqref{eqn:ortho_proj} and \eqref{eqn:sum_proj} that
		\begin{gather*}
			\PP_{\flat, j}(\eps D) = \PP_{\flat, j}^{(0)} \left( \widetilde{D} \right) + \eps |D| \PP_{\flat, j}^{(1)}(\eps D),\\
			\PP_{\flat, j}^{(0)} \Pi = \Pi \PP_{\flat, j}^{(0)}= \PP_{\flat, j}^{(0)},
		\end{gather*}
		therefore, there holds
		\begin{equation*}
			e^{t \lambda (\eps D) / \eps^2} \PP_{\flat, j}(\eps D) \left(\Id - \Pi\right)= e^{t \lambda (\eps D) / \eps^2} \eps |D| \PP_{\flat, j}^{(1)}(\eps D) \left(\Id - \Pi\right).
		\end{equation*}
		Thanks to \eqref{eqn:bound_hyproj_derivative_Linf}, we deduce the two estimates
		\begin{gather*}
			\left\|e^{t \lambda (\eps D) / \eps^2} \PP_{\flat, j}(\eps D) \left(\Id - \Pi\right) f \right\|_{\Spcg{\beta}} \lesssim  \left\| \chi(\eps D) \eps |D| e^{t \lambda (\eps D) / \eps^2} \right\|_{\BBB({\Spcg{\beta}})} \|f\|_{\Spcp{p}{\alpha}},\\
			\left\|e^{t \lambda (\eps D) / \eps^2} \PP_{\flat, j}(\eps D) \left(\Id - \Pi\right) f \right\|_{\Spcg{\beta}} \lesssim \left\| \eps \chi(\eps D) e^{t \lambda (\eps D) / \eps^2} \right\|_{\BBB({\Spcg{\beta}})} \|\nabla_x f\|_{\Spcp{p}{\alpha}}.
		\end{gather*}
		The operators $\chi(\eps D) \eps |D| e^{t \lambda (\eps D) / \eps^2}$ and $\chi(\eps D) e^{t \lambda (\eps D) / \eps^2} $ are Fourier multiplier in $x$, and \eqref{eqn:eigen_expansion} implies that $\Re \lambda_j(\xi) \leq - c |\xi|^2$ for some $c > 0$, thus, recalling that $\chi$ is the characteristic function of $\{ |\xi| \leq \err \}$, we have on the one hand
		\begin{align*}
			\left\| \chi(\eps D) \eps  |D| e^{t \lambda (\eps D) / \eps^2} \right\|_{\BBB({\Spcg{\beta}})} &\lesssim \sup_{|\eps \xi| \leq \err} \left|\eps \xi e^{-c t |\xi|^2}\right|\\
			& \lesssim \min\left\{ \frac{\eps}{t^{1/2}}, 1 \right\},
		\end{align*}
		and on the other hand
		$$\left\| \chi(\eps D) e^{t \lambda (\eps D) / \eps^2} \right\|_{\BBB({\Spcg{\beta}})}  \lesssim 1.$$
		These estimates yield the conclusion thanks to \eqref{eqn:semigroup_branch_def}.
	\end{proof}
	
	\begin{lem}
		\label{lem:cross_partial_semigroup}
		For any $p\in [1, \infty]$, $\alpha > \alpha_\BB$, $\beta \geq 0$, there holds for any $f \in \Spcp{p}{\alpha}$
		\begin{gather}
			\label{eqn:vanish_semigroup_hydro}
			\lim\limits_{\eps \to 0} \sup_{t \geq 0} \bigg(\langle t \rangle^{1/2} \|U^\eps_\flat(t) \left(\textnormal{Id} - \Pi\right) f\|_{\Spcg{\beta}}\bigg) = 0,\\
			\label{eqn:vanish_semigroup_hydro_ortho}
			\lim\limits_{\eps \to 0} \sup_{t \geq 0} \bigg(\langle t \rangle^{1/2} \|U^\eps_\sharp(t) \Pi f\|_{\Spcg{\beta}} \bigg)= 0.
		\end{gather}
	\end{lem}

	\begin{proof}
		First, note that $\alpha_\BB(1) = 2$ and $\alpha_\BB(\infty) > 5$, thus, when $p=\infty$, there is some~$\alpha' > 2$ such that $\Spcp{\infty}{\alpha} \subset \Spcp{1}{\alpha'}$ continuously, so we assume $p < \infty$. Second, since Lemma \ref{lem:est_u_pi_perp} implies that
		\begin{equation*}
			\langle t \rangle^{1/2} \|U^\eps_\flat(t) \left(\textnormal{Id} - \Pi\right) f\|_{\Spcg{\beta}} \lesssim \|f\|_\Spcp{p}{\alpha},
		\end{equation*}
		it is enough to check that the convergence holds on a dense subset of $\Spcp{p}{\alpha}$, for instance $\CC^\infty_c \left( \R_v \times \Omega_x \right)$. This is indeed the case for \eqref{eqn:vanish_semigroup_hydro} thanks to Lemma \ref{lem:est_u_pi_perp}.
		
		\smallskip
		The proof of \eqref{eqn:vanish_semigroup_hydro_ortho} is similar, one just needs to notice that
		\begin{align*}
			\langle t \rangle^{1/2} \|U^\eps_\sharp(t) \Pi f\|_{\Spcg{\beta}} &=  \langle t \rangle^{1/2} \|U^\eps_\sharp(t) \hyprojo \Pi f\|_{\Spcg{\beta}}\\
			& \lesssim e^{-\ah t / \eps^2} \langle t \rangle^{1/2} \|\hyprojo \Pi f\|_{\Spcg{\beta}} \\
			& \lesssim \|\left(\hyproj - \Pi \right) \Pi f\|_{\Spcg{\beta}}.
		\end{align*}
		Moreover, we have from \eqref{eqn:proj_expansion} and \eqref{eqn:sum_proj}
		\begin{equation*}
			\hyproj - \Pi = \sum_{j=-1}^{2} \eps |D| \PP_{\flat, j}^{(1)}(\eps D),
		\end{equation*}
		so by a similar density argument as in the previous step, one shows that $\left(\hyproj - \Pi \right) g$ vanishes as $\eps$ goes to zero for any~$g \in \Spcp{p}{\alpha}$, which concludes the proof.
		
	\end{proof}
	
	\section{Study in the Gaussian space}
	\resetcounter
	\label{scn:a_priori_gauss}
	Let us define the threshold appearing in Theorem \ref{thm:main} as
	\begin{equation}
		\label{eqn:def_alpha_star}
		\alpha_*(p) := \max \{\alpha_Q, \alpha_\BB\} + 1,
	\end{equation}
	where $\alpha_Q, \alpha_\BB$ are defined respectively in \eqref{eqn:def_alpha_Q} and \eqref{eqn:def_alpha_B}. In particular, $\alpha_*(1) = 3$.
	\medskip
	
	In this section, we prove the estimates necessary to the study of the equation on~$g^\eps$ in \eqref{eqn:coupled_diff}. Namely, we show that $g \mapsto \Psi^\eps(f, g)$ (Lemma \ref{lem:bound_psi}-\ref{lem:bound_psi_NSFI}) and $\Phi^\eps[h]$ (Lemma~\ref{lem:est_phi}) have small operator norms in several cases (depending on $f$). Following the idea of \cite{IGT20}, the particular case of $f=f^0$ being the solution to \eqref{eqn:KINSF} will be dealt with by introducing an equivalent norm, which we introduce now. Several results are already present in \cite{IGT20}, we prove them for the sake of completeness.
	\begin{notation}
		\label{not:func_spcg}
		For any $\lambda > 0$ and $\beta \geq 0$, we denote $\tspcg{\lambda}{\beta}$ the set of continuous functions $g \in \CC_b\left( [0, T) ; \Spcg{\beta} \right)$ satisfying $\tspcgN{g}{\lambda}{\beta} < \infty$, where
		\begin{gather*}
			\tspcgN{g}{\lambda}{\beta} := \sup_{0 \leq t < T} \left\|\Lambda(t, \lambda) \chi_\Omega(t) g(t) \right\|_{\Spcg{\beta}},\\
			\Lambda(t, \lambda) := \exp\left( - \lambda \int_0^t \|f^0(t')\|_{_{\Spcg{\beta}}}^3 \d t'\right),
		\end{gather*}
		and we denoted $\chi_\Omega(t) := \langle t \rangle^{1/4}$ if $\Omega = \R^2$, and $\chi_\Omega(t) := 1$ otherwise.
	\end{notation}

	The function $t \mapsto \|f^0(t)\|_\Spcg{\beta}$ lies in $L^\infty \cap L^2$ according to Theorem \ref{thm:existence_INSF} and thus, in particular, in $L^3$ (any $q > 2$ would suffice but we chose the smallest admissible integer for clarity, see Remark \ref{rem:why_3}). This norm therefore satisfies for some $C=C(\lambda) > 0$ and any $f \in \tspcg{\lambda}{\beta}$
	\begin{equation}
		\label{eqn:equi_norm}
		\|f\|_{L^\infty_t \Spcg{\beta} } \leq \|\chi_\Omega f\|_{L^\infty_t \Spcg{\beta}} \leq \frac{1}{C} \tspcgN{f}{\lambda}{\beta} \leq C \|\chi_\Omega f\|_{L^\infty_t \Spcg{\beta}},
	\end{equation}
	showing that $f \mapsto \|\chi_\Omega f \|_{L^\infty_t \Spcg{\beta}}$ and $f \mapsto \tspcgN{f}{\lambda}{\beta}$ are equivalent norms. The factor $\Lambda$ allows not to assume $f^0$ (and thus $f_{\ini, \WP}$) to be small thanks to the relation
	\begin{equation}
		\label{eqn:int_NSFI}
		\Lambda(t, \lambda)^3 \int_{0}^{t} \left\| \Lambda(\lambda, t')^{-1} f^0(t')\right\|^3_\Spcg{\beta} \d t' \lesssim \lambda^{-1},
	\end{equation}
	and the factor $\chi_\Omega$ makes the bilinear operator $\Psi^\eps$ bounded.
	
	\begin{lem}
		\label{lem:bound_psi}
		For any $\lambda > 0$, $\beta > d/2 + 1$ and $\omega > 0$, there holds for some positive constant $C=C(\lambda, \omega)$
		\begin{gather}
			\label{eqn:bound_psi_bil}
			\tspcgN{\Psi^\eps(f, g)}{\lambda}{\beta} \leq C \tspcgN{f}{\lambda}{\beta} \tspcgN{g}{\lambda},\\
			\label{eqn:bound_psi_exp}
			\left\| \Psi^\eps (f, g) \right\|_{\tspcg{\lambda}{\beta}} \leq C \eps \tspcgN{f}{\lambda}{\beta} \bigg(\sup_{0 \leq t < T} \left\|e^{\omega t / \eps^2} g(t)\right \|_{\Spcg{\beta}}\bigg).
		\end{gather}
	\end{lem}

	\begin{proof}
		
		\step{1}{Reductions}
		Using the same notations, a factorization similar to the one used in the proof of \eqref{eqn:Linf_dec_SG_hydroo} holds:
		\begin{align}
			\notag
			\Psi^\eps(f, g) &= \frac{1}{\eps} U^\eps * Q(f, g) \\
			\notag
			&= \frac{1}{\eps} S^\eps * Q(f, g) + \frac{1}{\eps^2} S^\eps * KU^\eps * Q(f, g) \\
			\label{eqn:fact_psi}
			&=: \Psi_0^\eps(f, g) + \frac{1}{\eps^2} S^\eps * K \Psi^\eps(f, g).
		\end{align}
		Again, denoting $(X, Y) = (\spcg^\gamma , \spcg^{\gamma+1})$ for any~$\gamma \geq 0$, or $(X, Y)=( H, \spcg^0)$ (where we recall $H=H^s_x L^2_v \left(M^{-1/2} \right)$), the regularization property $K \in \BBB\left(X; Y\right)$ and the decay of $S^\eps$ imply that:
		\begin{align*}
			\left\|\chi_\Omega S^\eps * K \Psi^\eps(t)(f, g) \right\|_{Y} &\lesssim \left( \chi_\Omega(t) \int_0^t \frac{e^{-\nu_0 (t-t')/\eps^2}}{\chi_\Omega(t')} \d t' \right) \|\chi_\Omega \Psi^\eps(f, g)\|_{L^\infty_t X} \\
			& \lesssim \eps^2 \|\chi_\Omega \Psi^\eps(f, g)\|_{L^\infty_t X},
		\end{align*}
		where we used Lemma \ref{lem:int_est_no_sing} in the second line. By induction, it is then enough to prove
		\begin{equation}
			\label{eqn:fact_psi_norm}
			\|\chi_\Omega \Psi^\eps(f, g)\|_{L^\infty_t \Spcg{\beta}} \lesssim \|\chi_\Omega  \Psi_0^\eps(f, g)\|_{L^\infty_t \Spcg{\beta}} +  \|\chi_\Omega \Psi^\eps(f, g)\|_{L^\infty_t \spcgg}.
		\end{equation}
	\step{2}{Estimate of $\Psi_0^\eps$}
	Note that for any $h=h(x, v)$, by \eqref{eqn:semigroup_transport}, one has
	$$\|S^\eps(t) h(v)\|_{H_x^s} = \exp(-\nu(v) t / \eps^2) \|h\|_{H^s_x},$$
	thus, we have
	\begin{align*}
		\chi_\Omega(t)\langle v \rangle^{\beta} \|\Psi_0^\eps(f, g)(t, v)\|_{H^s_x} &\leq \frac{1}{\eps} \langle v \rangle^{\beta} \chi_\Omega(t) \int_{0}^{t} e^{ -\nu(v) (t-t') / \eps^2 } \|Q(f, g)(t', v)\|_{H^s_x} \d t'\\
		&\leq \frac{1}{\eps} \langle v \rangle \chi_\Omega(t) \int_{0}^{t} e^{ -\nu(v) (t-t') / \eps^2 } \|f(t')\|_{\Spcg{\beta}} \|g(t')\|_{\Spcg{\beta}} \d t'\\
		& \leq \left( \frac{1}{\eps} \chi_\Omega(t) \int_{0}^{t} \frac{e^{ -\nu(v) (t-t') / \eps^2 } \langle v \rangle }{\chi_\Omega(t')^2} \d t' \right) \|\chi_\Omega f\|_{L^\infty_t \Spcg{\beta}} \|\chi_\Omega g\|_{L^\infty_t \Spcg{\beta}}
	\end{align*}
	where we used in the second line the following estimate from \cite[(B.5)]{IGT20}
	\begin{equation*}
		\|Q(f, g)\|_{\spcg^\gamma} \lesssim  \|f\|_{\spcg^{\gamma+1}} \|g\|_{\spcg^{\gamma+1}}.
	\end{equation*}
	The factor between parenthesis is bounded uniformly in $t$, $\eps$ and $v$ by Lemma \ref{lem:int_est_no_sing} and~\eqref{eqn:bound_col_freq}, therefore we obtain
	\begin{equation}
		\label{eqn:bound_phi}
		\|\chi_\Omega \Psi_0^\eps(f, g)\|_{L^\infty_t \Spcg{\beta}} \lesssim \eps \| \chi_\Omega f \|_{L^\infty_t \Spcg{\beta}} \| \chi_\Omega g \|_{L^\infty_t \Spcg{\beta}}.
	\end{equation}

	\step{2}{Proof of \eqref{eqn:bound_psi_bil}}
	Recall that the microscopic laws of elastic collisions imply~$\Pi Q = 0$ (see for instance \cite[(1.2.7)]{UY06}), thus one may write
	\begin{gather*}
		\Psi^\eps(f, g) = \frac{1}{\eps} U^\eps * Q(f, g) = W^\eps * Q(f, g),\\
		W^\eps := \frac{1}{\eps} U^\eps \left(\Id - \Pi \right).
	\end{gather*}
	By using the decay estimate \cite[Lemma 3.2]{IGT20} of $W^\eps$ then the bounds \cite[(B.5)-(B.6)]{IGT20} on $Q(f, g)$, one has
	\begin{align*}
		\|W^\eps(t-t') Q(f(t'), g(t'))\|_\spcgg & \lesssim \widetilde{\chi}_\Omega(t-t')\|Q(f(t'), g(t'))\|_{ \spcgg \cap L^2_v L^1_x \left( M^{-1/2} \right) }\\
		& \lesssim \widetilde{\chi}_\Omega(t-t') \|f(t')\|_{\Spcg{\beta}} \|g(t')\|_{\Spcg{\beta}}.
	\end{align*}
	where we have denoted for some $\sigma > 0$
	\begin{equation}
		\label{eqn:chi_omega_tilde}
		\widetilde{\chi}_\Omega(t) := 
		\begin{cases}
			t^{-1/2} e^{-\sigma t}, &  \Omega = \T^d,\\	 
			t^{-1/2} \langle t \rangle^{-d/4}, & \Omega = \R^d,\\	
		\end{cases}
	\end{equation}
	We thus have the control
	\begin{align*}
		\left\|\chi_\Omega \Psi^\eps(f, g)\right\|_{L^\infty_t \spcgg} & \lesssim \int_0^t \|\chi_\Omega(t) W^\eps(t-t') Q(f(t'), g(t'))\|_\spcgg \, \d t'\\
		&\lesssim \left( \chi_\Omega(t)\int_0^t \frac{ \widetilde{\chi}_\Omega(t-t') }{ \chi_\Omega(t')^2} \d t' \right) \|\chi_\Omega f\|_{L^\infty_t {\Spcg{\beta}}} \|\chi_\Omega g\|_{L^\infty_t {\Spcg{\beta}}}\\
		& \lesssim \|\chi_\Omega f\|_{L^\infty_t {\Spcg{\beta}}} \|\chi_\Omega g\|_{L^\infty_t {\Spcg{\beta}}},
	\end{align*}
	where the last line comes from Lemma \ref{lem:int_est_general}. Combined with \eqref{eqn:bound_phi}, this yields~\eqref{eqn:fact_psi_norm} hence \eqref{eqn:bound_psi_bil}.
	
	\step{3}{Proof of \eqref{eqn:bound_psi_exp}}
	Similarly, we have
	\begin{align*}
		\|\chi_\Omega \Psi^\eps(f, g)\|_{L^\infty_t \spcgg} & \lesssim \int_0^t \|\chi_\Omega(t) W^\eps(t-t') Q(f(t'), g(t'))\|_\spcgg \, \d t'\\
		&\lesssim \left( \chi_\Omega(t) \int_0^t \frac{ \widetilde{\chi}_\Omega(t-t') e^{-\omega t' / \eps^2} }{ \chi_\Omega(t')} \right) \|\chi_\Omega f\|_{L^\infty_t {\Spcg{\beta}}} \bigg(\sup_{0 \leq t' \leq T} \left\|e^{\omega t' / \eps^2} g(t')\right \|_{{\Spcg{\beta}}}\bigg).
	\end{align*}
	Thanks to Lemma \ref{lem:int_est_general} and the control $e^{-\omega t / \eps^2} \lesssim \eps t^{-1/2}$, this implies
	\begin{equation*}
		\|\chi_\Omega \Psi^\eps(f, g)\|_{L^\infty_t \spcgg} \lesssim \eps \|\chi_\Omega f\|_{L^\infty_t {\Spcg{\beta}}} \bigg(\sup_{0 \leq t' \leq T} \left\|e^{\omega t' / \eps^2} g(t')\right \|_{{\Spcg{\beta}}}\bigg)
	\end{equation*}
	Again, combined with \eqref{eqn:fact_psi_norm} and \eqref{eqn:bound_phi}, we deduce \eqref{eqn:bound_psi_exp}.

%

	\end{proof}
	
	The next lemma deals with the part of the source term $\SS^\eps$ from Section \ref{scn:main_reductions} associated with the acoustic waves. These waves, generated by the ill-prepared part~$f_{\ini, \IP}$ of the initial data, satisfy dispersive estimates when $\Omega=\R^d$ (see \eqref{eqn:Linf_statphase_disp_bis}), and are absent when $\Omega=\T^d$ thanks to the well-prepared assumption $f_{\ini, \IP} = 0$.
	
	\begin{lem}
		\label{lem:vanish_psi_disp}
		Under the assumptions of Theorem \ref{thm:main}, for any $\lambda > 0$, and $\beta > d/2+1$, there holds
		\begin{gather*}
			\begin{matrix}
				\left\|\Psi^\eps \left( U^\eps_{\disp} f_{\ini}, \, \cdot \, \right)\right\|_{\BBB\left(\tspcg{\lambda}{\beta}\right)} \xrightarrow[\eps \to 0]{} 0, &~ \text{if } \Omega = \R^d,\\
				\Psi^\eps \left( U^\eps_{\disp} f_{\ini}, \, \cdot \, \right) = 0, &~ \text{if } \Omega = \T^d.
			\end{matrix}
		\end{gather*}
	\end{lem}

	\begin{proof}
		Note that $U^\eps_\disp f_\ini = U^\eps_\disp f_{\ini, \IP}$ by \eqref{eqn:proj_hydro_semigroup}, so  the well-prepared assumption in the case $\Omega = \T^d$ concludes the proof in this case. We only deal with case $\Omega = \R^d$ in the following.

		\step{1}{Reductions}
		Let us denote by $\mathbf{X}$ the space
		\begin{equation*}
			\mathbf{X} := 
			\begin{cases}
				H^s_x, & \text{ if } \Omega = \R^3,\\
				H^s_x \cap L^1_x, & \text{ if } \Omega = \R^2,
			\end{cases}
 		\end{equation*}
		As $U^\eps_\disp f_\ini = U^\eps_\disp f_{\ini, \IP}$ by \eqref{eqn:proj_hydro_semigroup} and $f_{\ini, \IP}$ is macroscopic in the sense of Notations~\ref{scn:micro_macro_decomposition}, it writes
		$$f_{\ini}(x, v)= \bigg(\rho(x) + u(x) \cdot v + \frac{\theta(x)}{2} (|v|^2 - d)  \bigg)M(v)$$
		for some functions $\rho, u, \theta \in \mathbf{X}$. Consider a sequence of $\CC^\infty_c$ functions $\rho_n, u_n, \theta_n$ converging to $\rho, u, \theta$ in $\mathbf{X}$, and denote $f_n$ the corresponding macroscopic distribution. In both cases, there holds thanks to \eqref{eqn:Linf_bound_disp} and~\eqref{eqn:Linf_decay_disp}
		\begin{equation*}
			\chi_\Omega(t) \tspcgN{U^\eps_\disp (f_{\ini, \IP} - f_n)}{\lambda}{\beta} \lesssim \|f_{\ini, \IP} - f_n\|_{L^\infty_v \mathbf{X}_x \left( \langle v \rangle^\beta M^{-1/2} \right) } \xrightarrow[n \to \infty]{} 0.
		\end{equation*}
		Therefore, by the continuity \eqref{eqn:bound_psi_bil} of $\Psi^\eps$ on $\tspcg{\lambda}{\beta}$, one only need to check that for each~$n \in \N$,
		\begin{equation*}
			\left\|\Psi^\eps\left(U^\eps_\disp f_n, \, \cdot \, \right)\right\|_{\BBB\left(\tspcg{\lambda}{\beta}\right)} \xrightarrow[\eps \to 0]{} 0.
		\end{equation*}
		Note that by \eqref{eqn:Linf_bound_disp}, there holds~$U^\eps_\disp f_{\ini, \IP} \in \tspcg{\lambda}{\beta}$ thus we do not need to deal with~$\Psi_0^\eps$ as we have already proved \eqref{eqn:bound_phi}.

		\step{2}{Convergence for smooth functions}
		Thanks to the following bound for $Q$ from \cite[(B.4)]{IGT20}:
		\begin{equation*}
			\|Q(f, g)\|_{H^s_x L^2_v \left( M^{-1/2} \right)} \lesssim  \|f\|_{L^\infty_v W^{s, \infty}_x \left( \langle v \rangle^\beta M^{-1/2} \right) } \|g\|_{\Spcg{\beta}}
		\end{equation*}
		we can use \eqref{eqn:Linf_statphase_disp_bis}, thus we have
		\begin{equation*}
			\tspcgN{\Psi_0^\eps \left( U^\eps_\disp f_n, g \right)}{\lambda}{\beta} \lesssim C_{f_n} \eps^{1/2} \left( \chi_\Omega(t) \int_{0}^{t} \frac{\widetilde{\chi}_\Omega(t-t') }{\chi_\Omega(t') (t')^{1/2}} \d t' \right) \tspcgN{g}{\lambda}{\beta},
		\end{equation*}
		where $C_{f_n}$ is the constant of \eqref{eqn:Linf_statphase_disp_bis} and $\widetilde{\chi}_\Omega$ was defined in \eqref{eqn:chi_omega_tilde}. The integral is bounded uniformly in $t$ by Lemma \ref{lem:int_est_general}. We conclude the proof thanks to \eqref{eqn:fact_psi_norm} and~\eqref{eqn:bound_phi}.
	\end{proof}
	
	The next lemma shows how the equivalent norm $\tspcgN{\cdot}{\lambda}{\beta}$ defined in Notations \ref{not:func_spcg} makes the norm of the operator $\Psi^\eps(f^0, \, \cdot \,)$ arbitrarily small when $\eps \ll 1$ and $\lambda \gg 1$. Thisis why we do not need to assume that the well-prepared part of the initial data (which generates $f^0$) to be small.
	
	\begin{lem}
		\label{lem:bound_psi_NSFI}
		For any $\beta > d/2+1$, there exists a constant $C > 0$ (independent of~$\lambda$ and $\eps$) such that
		\begin{equation*}
			\tspcgN{\Psi^\eps(f^0, g)}{\lambda}{\beta} \leq C \left(\lambda^{-1/5} + \eps\right) \tspcgN{g}{\lambda}{\beta}.
		\end{equation*}
	\end{lem}

	\newcommand{\tspcgg}[1]{Y^{#1}}
	\newcommand{\tspcgNg}[2]{\left\|#1\right\|_{\tspcgg{#2}} }

	\begin{proof}
		Note that $t \mapsto \Lambda(\lambda, t)$, defined in Notations \ref{not:func_spcg}, is a decreasing function, thus for any positive functions $\phi=\phi(t)$, $\psi=\psi(t)$, one has
		\begin{equation}
			\label{eqn:convo_trick}
			\Lambda(\lambda, \cdot) \bigg(\phi * \psi\bigg) \leq \phi * \bigg( \Lambda(\lambda, \cdot) \psi \bigg),
		\end{equation}
		which implies that the bootstrap leading to \eqref{eqn:fact_psi_norm} is still valid for the norm of $\tspcg{\lambda}{\beta}$, uniformly in $\lambda$:
		\begin{align}
			\tspcgN{\Psi^\eps(f, g)}{\lambda}{\beta} &\lesssim \tspcgN{\Psi_0^\eps(f, g)}{\lambda}{\beta} + \|\chi_\Omega \Lambda(\lambda, \cdot) \Psi^\eps(f, g)\|_{L^\infty_t \spcgg }.
		\end{align}
		The term involving $\Psi_0^\eps$ is estimated the same way as for \eqref{eqn:bound_phi} thanks to \eqref{eqn:convo_trick}:
		\begin{equation*}
			\tspcgN{\Psi_0^\eps(f^0, g)}{\lambda}{\beta} \lesssim \eps \tspcgN{g}{\lambda}{\beta}.
		\end{equation*}
		The same goes for $\Psi^\eps$; the proof of \eqref{eqn:bound_psi_bil} leads to the bound
		\begin{align*}
			\tspcgN{\Psi^\eps(f^0, g) }{\lambda}{\beta}&\lesssim \sup_{t \geq 0} \left( \Lambda(\lambda, t) \chi_\Omega(t) \int_{0}^{t} \widetilde{\chi}_\Omega(t-t')  \|f^0(t')\|_\Spcg{\beta} \|g(t')\|_{\Spcg{\beta}} \, \d t' \right)\\
			& \lesssim \sup_{t \geq 0} \left( \int_{0}^{t} \frac{\chi_\Omega(t) \widetilde{\chi}_\Omega(t-t') }{\chi_\Omega(t')} \times \chi_\Omega(t') \Lambda(\lambda, t') \|f^0(t')\|_{\Spcg{\beta}} \, \d t' \right) \tspcgN{g}{\lambda}{\beta}
		\end{align*}
		where $\widetilde{\chi}_\Omega$ is defined in \eqref{eqn:chi_omega_tilde} and we used the fact that $t \mapsto \Lambda(\lambda, t)$ is non-increasing. The Hölder inequality with exponents $\left(3/2, 3\right)$ yields thanks to \eqref{eqn:int_NSFI}
		\begin{equation*}
			\tspcgN{\Psi^\eps(f^0, g) }{\lambda}{\beta} \lesssim \lambda^{-1/3} \tspcgN{g}{\lambda}{\beta},
		\end{equation*}
		which concludes the proof.
	\end{proof}

	\begin{rem}
		\label{rem:why_3}
		Note that in the previous proof, we use Hölder's inequality with exponents $(3/2, 3)$, but we could actually use any $(q, q')$ as long as $q \in (1, 2)$ so that the singularity of $\widetilde{\chi}_\Omega$ be integrable.
	\end{rem}

	Let us now study the convolution term $\eps^{-2} U^\eps * \AA$ appearing in the equation for~$g^\eps$ in \eqref{eqn:coupled_diff}. We start by introducing its splitting mentionned in Section \ref{scn:main_reductions}, then we prove estimates on both parts.
	\begin{notation}
		\label{not:convolution_split}
		For any $f=f(t, x, v)$, we denote
		\begin{gather*}
			\TT_1^\eps f := \frac{1}{\eps^2} U^\eps_\sharp * \AA f - \hyproj \exp\left( t \BB^\eps \right) f(0),\\
			\TT^\eps_\infty f := \frac{1}{\eps^2} U^\eps * \AA f - \TT^\eps_1 f,
		\end{gather*}
		where $U^\eps_\sharp$ and $\hyproj$ are defined in Section \ref{scn:spectral_splitting}. 
	\end{notation}

	\begin{lem}
		\label{lem:est_convo_1}
		Let $p \in [1, \infty]$, $\alpha > \alpha_\BB$ (defined in \eqref{eqn:def_alpha_B}), $\sigma \in (0, \sigma_\BB)$ (defined in Lemma \ref{lem:diss_B_eps}), $\beta > d/2$. There holds for some $C=C(p, \alpha, \beta, \lambda)$ and any $f \in \tSpcp{p}{\alpha}$
		\begin{equation}
			\label{eqn:est_convo_1}
			\| \TT^\eps_1 f (t)\|_{\Spcg{\beta}} \leq C e^{-\omega t/\eps^2 } \tSpcpN{f}{p}{\alpha},
		\end{equation}
		with $\omega := \min\{\sigma, \ah\}$, $\ah$ being defined in Section \ref{scn:spectral_splitting}.
	\end{lem}

	\begin{proof}
		Thanks to the decay estimate \eqref{eqn:Linf_dec_SG_hydroo} combined with the boundedness of teh operator $\AA : \Spcp{p}{\alpha} \rightarrow \Spcg{\beta}$, the definition of the norm $\tspcpN{f}$ from Notations~\ref{not:func_spc}, and the decay of $S_{\BB^\eps}$ coming from \eqref{eqn:coercivity}, we have
		\begin{align*}
			\|\TT^\eps_1 f (t)\|_{\Spcg{\beta}} &\lesssim \left\| \frac{1}{\eps^2} U^\eps_\sharp * \AA f \right\|_{\Spcg{\beta}} + \left\| \hyproj S_{\BB^\eps}(t) f(0) \right\|_{\Spcg{\beta}}  \\
			& \lesssim \frac{1 }{\eps^2} \tspcpN{f} \int_0^t e^{-\ah (t-t')/\eps^2} e^{-\sigma t'/\eps^2} \d t' + \|f(0)\|_{\Spcp{p}{\alpha}} \|S_{\BB^\eps}(t)\|_{\BBB(\Spcp{p}{\alpha})} \\
			& \lesssim e^{-\omega t / \eps^2} \tspcpN{f}.
		\end{align*}
		The lemma is proved.
	\end{proof}

	\begin{lem}
		\label{lem:est_convo_infty}
		Let $p \in [1, \infty]$, $\alpha > \alpha_*$ (defined in \eqref{eqn:def_alpha_star}), $\sigma \in (0, \sigma_\BB)$ (defined in Lemma \ref{lem:diss_B_eps}), $\beta > d/2 + 1$.
		For any $h \in L^\infty\left([0, T) ; \Spcp{p}{\alpha}\right) \cap L^p\left([0, T) ; \Spcpp{p}{\alpha}\right)$ and any~$g \in L^\infty([0, T) ; \Spcg{\beta})$, the solution~$\solp$ given by Lemma \ref{lem:a_priori_poly} to the equation
		\begin{equation}
			\label{eqn:eq_poly}
			\begin{cases}
				\partial_t\solp= \displaystyle \BB^\eps\solp+ \frac{1}{\eps} Q(h, h)  + \frac{1}{\eps} Q(h, g),\\
				\solp(0) \in \Spcp{p}{\alpha},
			\end{cases}
		\end{equation}
		is such that for some $C=C(p, \alpha, \beta, \lambda)$
		\begin{gather}
			\label{eqn:est_convo_infty}
			\tspcgN{\TT_\infty^\eps \solp}{\lambda}{\beta} \leq C \eps \tSpcpN{h}{p}{\alpha} \left(\tSpcpN{h}{p}{\alpha} + \tspcgN{g}{\lambda}{\beta} \right) + \tspcgN{U^\eps_\flat\solp(0)}{\lambda}{\beta}.
		\end{gather}
		Furthermore, if $\solp_1, \solp_2$ are two solutions of the equation \eqref{eqn:eq_poly} associated with $(h_1, g_1)$ and $(h_2, g_2)$ respectively, one has the stability estimate
		\begin{align}
			\notag
			\tspcgN{\TT_\infty^\eps \left(\solp_1 -\solp_2 \right)}{\lambda}{\beta} \leq& C \eps \tSpcpN{h_1 - h_2}{p}{\alpha} \left(\tSpcpN{h_1 + h_2}{p}{\alpha} + \tspcgN{g_1}{\lambda}{\beta} \right) \\
			\notag
			& + C \eps \tspcgN{g_1 - g_2}{\lambda}{\beta} \tSpcpN{h_2}{p}{\alpha} \\
			\label{eqn:est_convo_infty_stab}
			& + \tspcgN{U^\eps_\flat \left(\solp_1(0) -\solp_2(0)\right)}{\lambda}{\beta}.
		\end{align}
	\end{lem}
	
	\begin{proof}
		In the first step of the proof, which is entirely algebraic, we derive an expression of $\TT^\eps_\infty \solp$ in terms of $\solp(0)$, $h$ and $g$. In the next steps, we prove the estimates~\eqref{eqn:est_convo_infty} and \eqref{eqn:est_convo_infty_stab}. To lighten the notations, we denote
		\begin{equation*}
			\tspcp = \tSpcp{p}{\alpha}, ~ \spcp = \Spcp{p}{\alpha}, ~ \spcpp = \Spcpp{p}{\alpha},
		\end{equation*}
		where $\tSpcp{p}{\alpha}$ is defined in Notations \ref{not:func_spc}, and $\Spcp{p}{\alpha}, \Spcpp{p}{\alpha}$ in Notations \ref{not:func_spc_main}.
		We will also need the following factorization formula for semigroups \eqref{eqn:semigroup_fact} applied to the decomposition $\LL^\eps = \BB^\eps + \eps^{-2} \AA$: 
		\begin{gather}
			\label{eqn:duhamel_splitting_eps}
			U^\eps = S_{\BB^\eps} + \frac{1}{\eps^2} U^\eps * \AA S_{\BB^\eps},\\
			\notag
			S_{\BB^\eps}(t) := \exp\left(t \BB^\eps\right).
		\end{gather}
		
		\step{1}{Finding an expression for $\TT^\eps_\infty$}
		Using Duhamel's formula, $\solp$ writes
		\begin{equation*}
			\solp = S_{\BB^\eps}\solp(0) + \frac{1}{\eps} S_{\BB^\eps}* Q(h, h) + \frac{1}{\eps} S_{\BB^\eps}* Q\left(h, g\right),
		\end{equation*}
		and thus, convolving with $\eps^{-2} U^\eps \AA$, we have
		\begin{align*}
			\frac{1}{\eps^2} U^\eps * \AA\solp=& \frac{1}{\eps^2} U^\eps * \AA S_{\BB^\eps}\solp(0)\\
			&+ \frac{1}{\eps^2} U^\eps * \AA S_{\BB^\eps} * \left(\frac{1}{\eps} Q(h, h)\right) \\&+ \frac{1}{\eps^2} U^\eps * \AA S_{\BB^\eps} * \left(\frac{1}{\eps} Q(h, g)\right).
		\end{align*}
		Thanks to \eqref{eqn:duhamel_splitting_eps}, this expression rewrites as
		\begin{align*}
			\frac{1}{\eps^2} U^\eps * \AA\solp=& \left( U^\eps - S_{\BB^\eps} \right)\solp(0)\\
			&+ \left( U^\eps - S_{\BB^\eps} \right) * \left(\frac{1}{\eps} Q(h, h)\right) \\&+ \left( U^\eps - S_{\BB^\eps} \right) * \left(\frac{1}{\eps} Q(h, g)\right)
		\end{align*}
		or, in a more compact way
		\begin{equation}
			\label{eqn:pre_split}
			\frac{1}{\eps^2} U^\eps * \AA\solp= U^\eps \solp(0) - S_{\BB^\eps}\solp(0) + A_1 + A_2,
		\end{equation}
		where we denoted
		\begin{gather*}
			A_1 := -\frac{1}{\eps}  S_{\BB^\eps} * Q(h, h) - \frac{1}{\eps}   S_{\BB^\eps} * Q(h, g),\\
			A_2 := \frac{1}{\eps} U^\eps * Q(h, h) + \frac{1}{\eps} U^\eps * Q (h, g).
		\end{gather*}
		Using the definitions $U^\eps_\flat = \hyproj U^\eps$, $\hyprojo = \Id - \hyproj$ and $U^\eps_\sharp = \hyprojo U^\eps$ from Section \ref{scn:spectral_splitting}, we deduce from \eqref{eqn:pre_split}
		\begin{align*}
			\frac{1}{\eps^2} U^\eps * \AA\solp &= \hyproj \left(\frac{1}{\eps^2} U^\eps * \AA\solp\right) + \hyprojo \left(\frac{1}{\eps^2} U^\eps * \AA\solp\right)\\
			&= U^\eps_\flat \solp(0) - \hyproj S_{\BB^\eps}\solp(0) + \hyproj A_1 + \hyproj A_2 + \frac{1}{\eps^2} U^\eps_\sharp * \AA\solp,
		\end{align*}
		which allows to identify $\TT^\eps_\infty$ as
		\begin{gather*}
			\TT^\eps_\infty\solp= \hyproj  A_1 + \hyproj  A_2 + U^\eps_\flat\solp(0).
		\end{gather*}

		\step{2}{Estimate of $\TT^\eps_\infty$}
		According to Lemma~\ref{lem:a_priori_poly}, there holds
		\begin{equation*}
			\tspcpN{A_1} \lesssim \eps \tspcpN{h} \left(\tspcpN{h} + \tspcgN{g}{\lambda}{\beta}\right),
		\end{equation*}
		therefore, by the boundedness of $\hyproj : \spcp \rightarrow \spcg$ and $\exp(-\sigma t / \eps^2)  \chi_\Omega(t)\lesssim 1$, one has
		\begin{equation}
			\label{eqn:ctrl_A1}
			\tspcgN{\hyproj A_1}{\lambda}{\beta} \lesssim \eps \tspcpN{h} \left(\tspcpN{h} + \tspcgN{g}{\lambda}{\beta}\right).
		\end{equation}
		To estimate the term $\hyproj A_2$, recall that the laws of elastic collisions imply $\Pi Q = 0$ (see for instance \cite[(1.2.7)]{UY06}), thus, by Lemma \ref{lem:est_u_pi_perp} (because $\alpha - 1 > \alpha_\BB$ by \eqref{eqn:def_alpha_star}):
		\begin{equation*}
			\frac{1}{\eps} \|\chi_\Omega(t) U^\eps_\flat * Q(h, h)(t) \|_{ \spcg } \leq \chi_\Omega(t) \int_{0}^{t} (t-t')^{-1/2} \left\|Q(h(t'), h(t')) \right\|_{\Spcp{p}{\alpha - 1 }} \d t'
		\end{equation*}
		and then, by \eqref{eqn:bound_Q} and because $\alpha - 1/p > \alpha_Q$ (by \eqref{eqn:def_alpha_star}),
		\begin{align*}
			 \frac{1}{\eps} \|\chi_\Omega(t) U^\eps_\flat * Q(h, h)(t) \|_{ \spcg } &\lesssim \chi_\Omega(t) \int_{0}^{t} (t-t')^{-1/2} \|h(t')\|_{\Spcp{p}{\alpha - 1/p}} \|h(t')\|_{\Spcpp{p}{\alpha-1/p}} \d t'\\
			 & \lesssim \chi_\Omega(t) \int_{0}^{t} (t-t')^{-1/2} \|h(t')\|_\spcp^2 \d t'.
		\end{align*}
		Using the estimates
		\begin{gather*}
			\|h(t)\|_{\spcp} \lesssim e^{-\sigma t / \eps^2} \tspcpN{h} \lesssim \eps t^{-1/2} \tspcpN{h},\\
			\chi_\Omega(t) \|h(t)\|_{\spcp} \lesssim \tspcpN{h},
		\end{gather*}
		we also have thanks to Lemma \ref{lem:int_est_general}
		\begin{align*}
			\frac{1}{\eps} \|\chi_\Omega(t) U^\eps_\flat * Q(h, h)(t) \|_{ \spcg } &\lesssim \eps \tspcpN{h}^2 \left( \chi_\Omega(t) \int_{0}^{t} \frac{\d t' }{(t-t')^{1/2} (t')^{1/2} \chi_\Omega(t')} \d t'\right)\\
			&\lesssim \eps \tspcpN{h}^2.
		\end{align*}
		By performing the same computation for the term involving $g$, one can show
		\begin{equation}
			\label{eqn:ctrl_A2}
			\tspcgN{\hyproj A_2}{\lambda}{\beta} \lesssim \eps \tspcpN{h}  \left(\tspcpN{h} + \tspcgN{g}{\lambda}{\beta}\right).
		\end{equation}
		We conclude to the estimate of $\TT^\eps_\infty$ thanks to \eqref{eqn:ctrl_A1}-\eqref{eqn:ctrl_A2}.
		
	\step{3}{Stability estimate}
	To obtain the stability estimates, notice that
	\begin{equation*}
		\partial_t (\solp_1 -\solp_2) = \BB^\eps \left(\solp_1 -\solp_2 \right) + \frac{1}{\eps} Q\left(h_1 - h_2, h_1 + h_2 + g_1\right) + \frac{2}{\eps} Q \left( h_2, g_1 - g_2 \right).
	\end{equation*}
	By adapting the previous steps, we get the result. Lemma \ref{lem:est_convo_infty} is proved.	
	\end{proof}

	This next lemma provides estimates for the operator $\Phi^\eps[h]$ which was defined in Section~\ref{scn:main_reductions} for any $h=h(t, x, v)$ as
	$$\Phi^\eps[h] := 2 \Psi^\eps\left(  f^0 + \TT^\eps_1 h + U^\eps_\disp f_{\ini}, \, \cdot \, \right).$$

	\begin{lem}
		\label{lem:est_phi}
		Consider $p \in [1, \infty]$, $\alpha > \alpha_*$, $\sigma \in (0, \sigma_\BB)$, $\beta > d/2 + 1$ and~$A > 0$. There holds for any~$g \in \tspcg{\lambda}{\beta}$, $h \in \tSpcp{p}{\alpha}$ satisfying $\tSpcpN{h}{p}{\alpha} \leq A$
		\begin{equation*}
			\tspcgN{\Phi^\eps[h] g}{\lambda}{\beta} \leq C_{\eps, \lambda, A} \tspcgN{g}{\lambda}{\beta},
		\end{equation*}
		and the following stability estimate
		\begin{equation*}
			\tspcgN{\Phi^\eps[h_1] g_1 - \Phi^\eps[h_2] g_2}{\lambda}{\beta} \leq C_{\eps, \lambda, A} \left(\tspcgN{g_1 - g_2}{\lambda}{\beta} + \tSpcpN{h_1 - h_2}{p}{\alpha}\right),
		\end{equation*}
		for some constant $C_{\eps, \lambda, A} > 0$ that vanishes as $\max\{\eps , 1/ \lambda\} \to 0$ with $A$ fixed..
	\end{lem}
	
	\begin{proof}
		The conclusion is a direct consequence of Lemma \ref{lem:bound_psi_NSFI} (for $f^0$),
		\eqref{eqn:bound_psi_exp} combined with~Lemma \ref{lem:est_convo_1} (for $\TT^\eps_1 h^\eps$), and Lemma \ref{lem:vanish_psi_disp} (for $U^\eps_\disp f_{\ini}$). For the same reason, the stability estimate holds if one notices that
		\begin{equation*}
			\Phi^\eps[h_1]g_1 - \Phi^\eps[h_2]g_2 = \Phi[h_1](g_1 - g_2) + 2 \Psi^\eps(\TT^\eps_1 (h_1 - h_2), g_2).
		\end{equation*}
		This proves the lemma.
	\end{proof}

	This next lemma provides estimates for $\SS^\eps$, which we recall is defined in Section~\ref{scn:main_reductions} as
	\begin{gather*}
		\SS^\eps\left[\solp\right] := \SS^\eps_0 + \Psi^\eps\left( \TT^\eps_1 \solp + U^\eps_\disp f_{\ini, \IP} , 2 f^0 + \TT^\eps_1 \solp + U^\eps_\disp f_{\ini, \IP} \right) + \TT^\eps_\infty \solp\\
		\SS^\eps_0 := \left(U^\eps - U^\eps_\disp - U^0\right) \left(f_{\ini} - f_{\ini, \bot}\right) + \left(\Psi^\eps - \Psi^0\right)\left(f^0, f^0\right).
	\end{gather*}

	\begin{lem}
		\label{lem:cv_source_NSFI}
		Let $p \in [1, \infty]$, $\alpha > \max\{\alpha_\BB, \alpha_Q\}$, $\sigma \in (0, \sigma_\BB)$, $\beta > d/2 + 1$, and any~$A > 0$. For any $h \in L^\infty\left([0, T) ; \Spcp{p}{\alpha}\right) \cap L^p\left([0, T) ; \Spcpp{p}{\alpha}\right)$, any $g \in L^\infty([0, T) ; \Spcg{\beta})$ such that~$\tspcgN{g}{\lambda}{\beta} + \tSpcpN{h}{p}{\alpha} \leq A$, the solution $\solp$ given by Lemma \ref{lem:a_priori_poly} to the equation
		\begin{equation*}
			\begin{cases}
				\partial_t \solp = \displaystyle \BB^\eps \solp + \frac{1}{\eps} Q(h, h) + \frac{1}{\eps} Q(h, g),\\
				\solp(0) = f_{\ini, \bot},
			\end{cases}
		\end{equation*}
		is such that
		\begin{equation}
			\label{eqn:vanish_source}
			\tspcgN{\SS^\eps\left[\solp\right]}{\lambda}{\beta} \leq \eta_\eps(\lambda, A, p, \alpha),
		\end{equation}
		where $\eta_\eps(\lambda, A) \xrightarrow[\eps \to 0]{} 0$ with $\lambda$ and $A$ fixed. Furthermore, one has the stability estimate
		\begin{equation}
			\label{eqn:stab_source}
			\tspcgN{\SS^\eps[\solp_1] - \SS^\eps[\solp_2]}{\lambda}{\beta} \leq C \eps A \left(\tSpcpN{h_1 - h_2}{p}{\alpha} + \tspcgN{g_1 - g_2}{\lambda}{\beta}\right),
		\end{equation}
		where $C = C(p, \alpha, \lambda)$.
	\end{lem}
	
	\begin{proof}
		\step{1}{Proof of \eqref{eqn:vanish_source}}
		One may rewrite $\SS^\eps_0$ as
		\begin{align*}
			\SS^\eps_0 =& \left(U^\eps - U^\eps_\disp - U^\eps_\sharp - U^0\right) \Pi f_{\ini}  \\
			&+ U^\eps_\sharp \Pi f_{\ini} \\
			&+ \left(\Psi^\eps - \Psi^0\right)\left(f^0, f^0\right).
		\end{align*}
		The first line vanishes in $\tspcg{\lambda}{\beta}$ in virtue of \cite[Lemma 3.5]{IGT20}, the second one by \eqref{eqn:vanish_semigroup_hydro_ortho}, and the third one by \cite[Lemma 4.1]{IGT20}.
		
		\medskip
		The term involving $\Psi^\eps$ expands as
		\begin{align*}
			\Psi^\eps( \TT^\eps_1 \solp + U^\eps_\disp f_{\ini, \IP} , 2 f^0 + &\TT^\eps_1 \solp + U^\eps_\disp f_{\ini, \IP} )  =\\
			& \Psi^\eps\left( \TT^\eps_1 \solp, 2 f^0 + \TT^\eps_1 \solp + U^\eps_\disp f_{\ini, \IP} \right)\\
			&+ \Psi^\eps\left( U^\eps_\disp f_{\ini, \IP} , 2 f^0 + \TT^\eps_1 \solp + U^\eps_\disp f_{\ini, \IP} \right).
		\end{align*}
		The first term is estimated using Lemma \ref{lem:est_convo_1} and \eqref{eqn:bound_psi_exp}, the second one using Lemma~\ref{lem:vanish_psi_disp}.
		
		\medskip
		Finally, the term $\TT_\infty^\eps \solp$ is estimated using \eqref{eqn:est_convo_infty} and \eqref{eqn:vanish_semigroup_hydro}. This proves \eqref{eqn:vanish_source}.
		
		\step{2}{Proof of \eqref{eqn:stab_source}}
		Note that
		\begin{align*}
			\SS^\eps[\solp_1] - \SS^\eps[\solp_2] =& \TT^\eps_\infty \left( \solp_1 - \solp_2 \right) \\
			&+ \Psi^\eps\left( \TT^\eps_1\left(\solp_1 - \solp_2\right) , \TT^\eps_1\left(\solp_1 + \solp_2\right) \right) \\
			& + 2 \Psi^\eps\left( \TT^\eps_1\left(\solp_1 - \solp_2\right) , U^\eps_\disp f_{\ini, \IP} \right).
		\end{align*}
		Using \eqref{eqn:est_convo_infty_stab}, and \eqref{eqn:bound_psi_exp} with Lemma \ref{lem:est_convo_1}, we get
		\begin{equation*}
			\tspcgN{\SS^\eps[\solp_1] - \SS^\eps[\solp_2] }{\lambda}{\beta} \lesssim \eps A \left(\tSpcpN{h_1 - h_2}{p}{\alpha} + \tspcgN{g_1 - g_2}{\lambda}{\beta}\right),
		\end{equation*}
		which concludes the proof of \eqref{eqn:stab_source}, hence of Lemma \ref{lem:cv_source_NSFI}.
		
	\end{proof}
	
	\section{Proof of Theorem \ref{thm:main}}

	\label{scn:cauchy_theory}
	\newcommand{\Pspc}{\XXX_{A, \lambda, \eps}^{p, \alpha, \sigma, \beta}}
	\newcommand{\pspc}{\XXX_{A, \lambda, \eps}}
	\newcommand{\pspcN}[1]{\left\|#1\right\|_{\pspc}}
	Let us fix $p \in [1, \alpha]$, $\beta > d/2+1$, and recall that $f_\ini, f^0$ and $T$ are fixed. We now prove Theorem \ref{thm:main} using Banach's fixed point theorem. To do so, we introduce the metric space in which we work:
	\begin{align*}
		\pspc=\Pspc := \bigg\{ (h, &g) \in \,\tSpcp{p}{\alpha} \times \tspcg{\lambda}{\beta} \, : \, \tSpcpN{h}{p}{\alpha} \leq 2 \|f_{\ini, \bot}\|_{\Spcp{p}{\alpha}}, ~ \tspcgN{g}{\lambda}{\beta} \leq A, \\
		& \exists \left(\widetilde{h}, \widetilde{g}\right) \in \tSpcp{p}{\alpha} \times \tspcg{\lambda}{\beta}, ~ \tSpcpN{\widetilde{h}}{p}{\alpha} \leq 2 \|f_{\ini, \bot}\|_{\Spcp{p}{\alpha}}, ~ \tspcgN{\widetilde{g}}{\lambda}{\beta} \leq A,\\
		&\partial_t h = \BB^\eps h + \frac{1}{\eps} Q\left(\widetilde{h}, \widetilde{h}\right) + \frac{1}{\eps} Q\left(\widetilde{h}, \widetilde{g}\right), ~ h(0) = f_{\ini, \bot}\bigg\},
	\end{align*}
	endowed with the norm
	\begin{equation*}
		\|(h, g)\|_{\pspc} := \tSpcpN{h}{p}{\alpha} + \tspcgN{g}{\lambda}{\beta}.
	\end{equation*}
	Note that by Lemma \ref{lem:a_priori_poly}, for $\eps$ small enough, this space is non-empty as it contains for instance $\left( e^{t \BB^\eps} f_{\ini, \bot}, 0\right)$.
	
	\medskip
	The map $\Xi$ on wich we will apply the fixed-point theorem is defined as
	\begin{gather*}
		\Xi \, : \, \pspc \rightarrow \pspc, ~ \Xi(h, g) := \left( \solp, \solg \right),\\
		\partial_t \solp = \displaystyle \BB^\eps \solp + \frac{1}{\eps} Q\left( h, h\right) + \frac{2}{\eps} Q \left(h, g+ f^0 + \TT^\eps_1 h + U^\eps_{\disp} f_{\ini}\right), ~ \solp(0) = f_{\ini, \bot},\\
		\solg = \SS^\eps[h] + \Psi^\eps(g, g) + \Phi^\eps[h] g,
	\end{gather*}
	where we recall the notations
	\begin{gather*}
		\SS^\eps[h] := \SS^\eps_0 + \Psi^\eps\left( \TT^\eps_1 h + U^\eps_\disp f_{\ini} , 2 f^0 + \TT^\eps_1 h + U^\eps_\disp f_{\ini} \right) + \TT^\eps_\infty h,\\
		\SS^\eps_0 := \left(U^\eps - U^\eps_\disp - U^0\right) \Pi f_{\ini} + \left(\Psi^\eps - \Psi^0\right)\left(f^0, f^0\right),\\
		\Phi^\eps[h] := 2 \Psi^\eps\left(  f^0 + \TT^\eps_1 h + U^\eps_\disp f_{\ini}, \, \cdot \, \right).
	\end{gather*}
	
	\begin{proof}[Proof of Theorem \ref{thm:main}]
		We prove in Step 1 that $\Xi$ is a well-defined contraction of $\pspc$ so that the system \eqref{eqn:coupled_diff} has a unique solution $(h, g)$ by Banach's fixed point theorem. As explained in Sections \ref{scn:main_reductions}, the functions $f := h + g$ will therefore be a solution of the scaled Boltzmann equation \eqref{eqn:bolt_scl}. In Step 2, we show that this solution is the unique solution to \eqref{eqn:boltz} satisfying \eqref{eqn:integ_solution}. In Step 3, we define the terms~$u^\eps_*$ and prove they vanish in appropriate topologies.

		\step{1}{$\Xi$ is a well-defined contraction}
		Let $(h_j, g_j) \in \pspc$ for $j=1, 2$, and denote their images by $\Xi$ by $\left( \solp_j, \solg_j \right) := \Xi(h_j, g_j)$. According to Lemma \ref{lem:a_priori_poly} and \eqref{eqn:equi_norm}, there holds for some $C(\lambda) >0$
		\begin{gather*}
			\tSpcpN{ \solp_j}{p}{\alpha} \leq C(\lambda) \eps (1+A) + \|f_{\ini, \bot}\|_{\Spcp{p}{\alpha}},\\
			\tSpcpN{ \solp_1 - \solp_2 }{p}{\alpha}  \lesssim C(\lambda) (1+A) \eps \pspcN{(h_1, g_1) -(h_2, g_2)}.
		\end{gather*}
		Furthermore, according to Lemma \ref{lem:est_phi}, \eqref{eqn:vanish_source} and \eqref{eqn:bound_psi_bil}, we may assume that $C(\lambda)$ also satisfies
		\begin{align}
			\notag
			\tspcgN{ \solg_j }{\lambda}{\beta} &\lesssim \tspcgN{\SS^\eps[h^\eps]}{\lambda}{\beta} + \tspcgN{\Psi^\eps(g_j, g_j)}{\lambda}{\beta} + \tspcgN{ \Phi^\eps[h_j]g_j}{\lambda}{\beta}\\
			\label{eqn:pre_vanish}
			&\lesssim \eta_\eps(\lambda, 2 \|f_{\ini, \bot}\|_{\Spcp{p}{\alpha}} + A) + \bigg(C(\lambda) A + C_{\lambda, \eps, 2 \|f_{\ini, \bot}\|_\Spcp{p}{\alpha} + A} \bigg) A
		\end{align}
		\begin{align*}
			\tspcgN{ \solg_1 - \solg_2}{\lambda}{\beta} \leq& \tspcgN{\SS^\eps[h_1] - \SS^\eps[h_2]}{\lambda}{\beta} + \tspcgN{\Psi^\eps\left( g_1 + g_2, g_1 - g_2 \right)}{\lambda}{\beta} \\
			&+ \tspcgN{ \Phi^\eps[h_1] g_1 - \Phi^\eps[h_2] g_2 }{\lambda}{\beta}\\
			\lesssim & \big(\eps C(\lambda)(1+A) + C(\lambda) A + C_{\eps, \lambda, 2 \|f_{\ini, \bot}\|_{\Spcp{p}{\alpha}}} \big) \pspcN{(h_1, g_1) -(h_2, g_2)}.
		\end{align*}
		First, set $A=1$. By Lemma \ref{lem:est_phi}, one must choose $\eps \ll 1$ and $\lambda \gg 1$ so that the values~$C_{\eps, \lambda, 2 \|f_{\ini, \bot}\|_{\Spcp{p}{\alpha}}}$ and $C_{\eps, \lambda, 2 \|f_{\ini, \bot}\|_{\Spcp{p}{\alpha}} + A}$ are small. Then, thanks to Lemma \ref{lem:cv_source_NSFI}, there holds up to a reduction of $A$ and $\eps$ that $\tSpcpN{ \solp_j}{p}{\alpha} \leq 2 \|f_\ini\|_{\Spcp{p}{\alpha}}$ and $\tspcgN{ \solg_j }{\lambda}{\beta} \leq A$ (the map $\Xi$ is well defined) and $\pspcN{(\solp_1, \solg_1) -(\solp_2, \solg_2)} \leq \frac{1}{2} \pspcN{(h_1, g_1) -(h_2, g_2)}$ (the map $\Xi$ is a contraction).
		
		We conclude that $\Xi$ has a unique fixed point $(h^\eps, g^\eps) \in \pspc$, and thus the system~\eqref{eqn:coupled_diff} has a unique solution $(h^\eps, g^\eps) \in \pspc$. Therefore, we have constructed a solution $f^\eps := h^\eps + g^\eps$ to the Boltzmann equation \eqref{eqn:bolt_scl}.
	
		\step{2}{Uniqueness of the solution}
		To show uniqueness, we may assume that the solution $f^\eps$ constructed previously satisfies the bound
		\begin{equation}
			\label{eqn:ass_solution_bound}
			\forall t \in [0, T), ~  \|f^\eps(t)\|_{\spcp} < \frac{\sigma_\BB}{2 C \eps }
		\end{equation}
		where $C$ is the constant of Lemma \ref{lem:estimate_Q}. Furthermore, by the definitions of $\alpha_*$ \eqref{eqn:def_alpha_star} and $\alpha_Q$ \eqref{eqn:def_alpha_Q}, one has
		\begin{equation*}
			\alpha > \alpha_*(p) \geq 3 + \frac{d}{p'},
		\end{equation*}
		thus there exists $\gamma > \alpha_*(1) = 3$ such that $\Spcp{p}{\alpha} \subset \Spcp{1}{\gamma}$ continuously. We therefore assume in this step that $p = 1$.
		Consider now another solution to \eqref{eqn:bolt_scl}
		\begin{equation*}
			f^\eps_1 \in \CC_b\left([0, T) ; \spcp\right) \cap L^p\left([0, T) ; \Spcpp{p}{\alpha}\right)
		\end{equation*}
		such that $f^\eps_1(0) = f_\ini$ and define 
		\begin{gather*}
			D := f^\eps - f^\eps_1, \\S := f^\eps + f^\eps_1.
		\end{gather*}
		These functions satisfy the following equation:
		\begin{equation*}
			\partial_t D = \LL^\eps D + \frac{1}{\eps}Q(D, S) = \BB^\eps D + \frac{1}{\eps^2} \AA D + \frac{1}{\eps}Q(D, S).
		\end{equation*}
		Similar calulations as in the proof of Lemma~\ref{lem:a_priori_poly} yield for some constant $C > 0$
		\begin{align*}
			\frac{\d}{\d t} \| D \|_{\spcp} \leq -\frac{\sigma_\BB}{\eps^2} \|D\|_{\spcpp} + \frac{C}{\eps} \left( \|D\|_{\spcp}  \|S\|_{\spcpp} + \|D\|_{\spcpp}  \|S\|_{\spcp} + \frac{1}{\eps} \|D\|_\spcp \right),
		\end{align*}
		which may be rewritten as
		\begin{gather}
			\label{eqn:diff_in}
			\frac{\d}{\d t} \bigg( e^{-\phi} \| D \|_{\spcp} \bigg) + \frac{ e^{-\phi} \|D\|_{\spcpp} }{\eps^2} \bigg( \sigma_\BB - C \eps \|S\|_{\spcp}\bigg) \leq 0,\\
			\notag
			\phi(t) := \frac{C}{\eps} \int_0^t \left\|S \left(t'\right) \right\|_{\spcpp} \d t' + \frac{C t}{\eps^2}.
		\end{gather}
		Define $I = \{ t \in [0, T) \, : \, D(t) = 0\}$. It is relatively closed in $[0, T)$ because $D$ is continuous, and non-empty since $0 \in I$. For any $t \in I$, assumption \eqref{eqn:ass_solution_bound} imply
		$$C \eps \|S(t)\|_{\spcp} = 2 C \eps \|f^\eps(t)\|_{\spcp} < \sigma_\BB,$$
		thus, by the continuity of $S$, there exists some $\delta > 0$ such that
		$$\forall t' \in [t, t+\delta), ~ C \eps \|S(t')\|_{\spcp} < \sigma_\BB.$$
		We conclude thanks to \eqref{eqn:diff_in} that $[0, t+\delta) \subset I$, thus $I$ is both relatively closed and open in $[0, T)$ and therefore $I = [0, T)$.
	
		\step{3}{Convergence of the terms $u_*^\eps$}
		We define the error terms of Theorem \ref{thm:main} in the following way:
		\begin{gather*}
			u^\eps_\infty(t) := g(t),\\
			u^\eps_1(t) := h(t) + \TT^\eps_1 h(t),\\
			u^\eps_\ac(t) := U^\eps_\disp(t) f_\ini.
		\end{gather*}
		The term $u^\eps_1$ satisfies by construction and Lemma \ref{lem:est_convo_1}
		\begin{equation*}
			\tSpcpN{u^\eps_1 }{p}{\alpha} \lesssim \|f_\ini\|_{\Spcp{p}{\alpha}},
		\end{equation*}
		which yields the estimate of Theorem \ref{thm:main} by the definition of the norm $\tSpcpN{\cdot}{p}{\alpha}$ in Notations \ref{not:func_spc}. The term $u^\eps_\ac$ vanishes in the senses stated in Theorem \ref{thm:main} thanks to the estimates of Lemma \ref{lem:disp}.
		Finally, $u^\eps_\infty$ vanishes uniformly because up to a reduction of the parameters $(G, 1/\lambda, \eps)$ in Step 1, one have from \eqref{eqn:pre_vanish} and Lemma \ref{lem:cv_source_NSFI}
		\begin{equation*}
			\tspcgN{u^\eps_\infty }{\lambda}{\beta} \lesssim \eta_\eps \rightarrow 0, ~ (\eps \to 0).
		\end{equation*}
	\end{proof}

	\appendix
	\section{Technical estimates}
	\resetcounter

		\begin{lem}
		\label{lem:int_est_no_sing}
		Let $a, \omega_0 > 0$, there exists some $C = C(a, \omega_0) > 0$ such that for any~$\omega \geq \omega_0$ and $t \geq 0$
		\begin{equation*}
			\int_0^t \frac{e^{- \omega t'}}{ \langle t-t' \rangle^a} \d t' \leq \frac{C}{\omega \langle t \rangle^a}.
		\end{equation*}
	\end{lem}
	
	\begin{proof}
		First, note that there holds
		\begin{equation*}
			\langle t \rangle \leq \langle t' \rangle \langle t-t' \rangle,
		\end{equation*}
		thus we have
		\begin{align*}
			\int_0^t \frac{e^{- \omega t'}}{ \langle t-t' \rangle^a} \d t' &\leq \langle t \rangle^{-a} \int_{0}^{t} e^{- \omega t'} \langle t' \rangle ^a \d t'\\
			& \lesssim \langle t \rangle^{-a} \int_{0}^{t} e^{- \omega t'} (1+ (t')^a) \d t'\\
			& \lesssim \frac{1}{\omega \langle t \rangle^a} + \langle t \rangle^{-a} \int_0^t e^{- \omega t'} (t')^a \d t'.
		\end{align*}
		Let us perform an integration by parts:
		\begin{align*}
			\int_0^t e^{- \omega t'} (t')^a \d t' &= \frac{a}{\omega} \int_0^t e^{- \omega t'} (t')^{ a-1 } \d t' - \frac{ e^{-\omega t}}{\omega} t^a\\
			& \lesssim \frac{1}{\omega} \left( \int_0^\infty (t')^{a-1} e^{-\omega_0 t'} \d t' - e^{- \omega_0 t}t^a \right).
		\end{align*}
		The factor between parenthesis is bounded because $a -1 > -1$ and thus the integral term converges.
	\end{proof}
	
	\begin{lem}
		\label{lem:int_est_general}
		For any $a, \alpha, c, \gamma \in [0, 1)$ and $b, \beta \geq 0$ such that 
		\begin{equation*}
			\begin{cases}
				a + \alpha \leq 1,\\
				c \leq a + b ,\\
				\gamma \leq \alpha + \beta,
			\end{cases}
		\end{equation*}
		there exists a constant $C > 0$ such that
		\begin{equation*}
			\int_{0}^t \frac{d t'}{ (t-t')^a \langle t - t' \rangle^b (t')^\alpha \langle t' \rangle^\beta } \leq C \langle t \rangle^{1 - c - \gamma}.
		\end{equation*}
	\end{lem}

	\begin{proof}
		Let us start with the change of variable $t = t' u$:
		\begin{align*}
			\int_{0}^t \frac{\d t'}{ (t-t')^a \langle t - t' \rangle^b (t')^\alpha \langle t' \rangle^\beta } & = t^{1 - a - \alpha} \int_{0}^1 \frac{\d u}{ (1-u)^a \langle t(1-u) \rangle^b u^\alpha \langle t u \rangle^\beta }.
		\end{align*}
		This quantity is bounded uniformly in $t \in [0, 1]$ because $a + \alpha \leq 1$. Furthermore,
		\begin{align*}
			\int_{0}^t \frac{\d t'}{ (t-t')^a \langle t - t' \rangle^b (t')^\alpha \langle t' \rangle^\beta } & \leq t^{1 - a - \alpha } \int_{0}^1 \frac{\d u}{ (1-u)^a \langle t(1-u) \rangle^{c-a} u^\alpha \langle t u \rangle^{\gamma - \alpha} },
		\end{align*}
		and note that for any $m \in \R$,
		\begin{equation*}
			t^m (t^{-1} + s)^m \lesssim \langle t s \rangle^m,
		\end{equation*}
		thus the integral can be controled by
		\begin{align*}
			\int_{0}^t \frac{\d t'}{ (t-t')^a \langle t - t' \rangle^b (t')^\alpha \langle t' \rangle^\beta } & \lesssim t^{1 - c - \gamma } \int_{0}^1 \frac{\d u}{ (1-u)^a (t^{-1} + 1-u)^{c-a}  u^\alpha (t^{-1} + u)^{\gamma-\alpha} }.
		\end{align*}
		The last thing to check is that this integral factor in the right hand side is bounded uniformly in $t \geq 1$; on the one hand
		\begin{equation*}
			(1-u)^{-a} (t^{-1} + 1-u)^{a-c} \leq
			\begin{cases}
				(1-u)^{-a} (2-u)^{a-c}, \text{ if } a \geq c,\\
				(1-u)^{-c}, \text{ if } a \leq c,
			\end{cases}
		\end{equation*}
		and on the other hand
		\begin{equation*}
			u^{-\alpha} (t^{-1} + u)^{\alpha - \gamma} \leq
			\begin{cases}
				u^{-\alpha} (1 + u)^{\alpha - \gamma}, \text{ if } \alpha \geq \gamma,\\
				u^{-\gamma}, \text{ if } \alpha \leq \gamma.
			\end{cases}
		\end{equation*}
		In each case, the integral converges because of the assumption $a, \alpha, c, \gamma \in [0, 1)$, which concludes the proof.
	\end{proof}
	
	\bibliographystyle{abbrv}
	\bibliography{bibliography}
	
\end{document}